\documentclass[a4paper,12pt]{article}
\usepackage{amssymb}
\usepackage{amsmath}
\usepackage{amscd}
\usepackage{amsthm}
\usepackage[arrow,matrix]{xy}
\usepackage[bookmarksnumbered,colorlinks,linkcolor=black,citecolor=black,urlcolor=black]{hyperref}
\usepackage{mathtools}

\usepackage{tikz-cd}
\usepackage[shortlabels]{enumitem}

\newtheorem{theo}{Theorem}[section]
\newtheorem{inttheo}{Theorem}
\newtheorem{prop}[theo]{Proposition}
\newtheorem{lem}[theo]{Lemma}
\newtheorem{cor}[theo]{Corollary}
\newtheorem{defi}[theo]{Definition}
\newtheorem{rem}[theo]{Remark}

\newtheorem{exe}[theo]{Exercise}
\newtheorem{exa}[theo]{Example}

\newtheorem{conj}[theo]{Conjecture}
\newtheorem{ques}[theo]{Question}

\theoremstyle{remark}

\newcommand{\bthe}{\begin{theo}}
\newcommand{\ble}{\begin{lem}}
\newcommand{\bpr}{\begin{prop}}
\newcommand{\bco}{\begin{cor}}
\newcommand{\bde}{\begin{defi}}
\newcommand{\ethe}{\end{theo}}
\newcommand{\ele}{\end{lem}}
\newcommand{\epr}{\end{prop}}
\newcommand{\eco}{\end{cor}}
\newcommand{\ede}{\end{defi}}
\newcommand{\brem}{\begin{rem}}
\newcommand{\erem}{\end{rem}}

\newcommand{\bexe}{\begin{exe}}
\newcommand{\eexe}{\end{exe}}

\newcommand{\bexa}{\begin{exa}}
\newcommand{\eexa}{\end{exa}}

\newcommand{\bconj}{\begin{conj}}
\newcommand{\econj}{\end{conj}}
\newcommand{\bques}{\begin{ques}}
\newcommand{\eques}{\end{ques}}

\def\Sel{{\rm Sel}}
\def\id{{\rm id}}

\def \qbar {{\overline{\Q}}}

\def \deg {{\rm{deg}}}
\def \Br {{\rm{Br}}}

\def \si {{\sigma}}
\def \Ga {{\Gamma}}
\def \R {{\mathbb{R}}}
\def \Pic {{\rm {Pic}}}
\def \Div {{\rm {Div}}}
\def \div {{\rm{div}}}
\def \Gal {{\rm{Gal}}}
\def \Ker {{\rm{Ker}}}
\def \Im {{\rm {Im}}}

\def \A{{\mathbb A}}
\def \P{{\mathbb P}}

\def \Spec {{\rm{Spec}}}
\def \dim {{\rm{dim}}}
\def \Hom {{\rm {Hom}}}
\def \End {{\rm {End}}}
\def \Pic {{\rm {Pic}}}
\def \GL {{\rm {GL}}}

\def \Aut{{\rm Aut}}

\def \Z {{\mathbb Z}}
\def \Q {{\mathbb Q}}
\def \F {{\mathbb F}}
\def \C {{\mathbb C}}
\def \Ext {{\rm Ext}}
\def \Id {{\rm Id}}

\def \val {{\rm{val}}}

\def \H {{\rm H}}
\def\G{{\mathbb G}}

\def\C{{\mathbb C}}

\def\lra{\longrightarrow}
\def\cl{{\rm cl}}

\def \s {{\rm s}}

\def\inv{{\rm inv}}

\def\NS{{\rm NS\,}}

\def\si{\sigma}
\def\cl{{\rm cl}}

\def\val{{\rm val}}

\def\Ga{\Gamma}

\def\et{{\rm{\acute et}}}

\newcommand{\RGamma}{\mathrm{R}\Gamma}
\newcommand{\Sym}{\mathrm{Sym}}

\newcommand{\Jac}{\mathrm{Jac}}
\newcommand{\proet}{\mathrm{proet}}
\newcommand{\cont}{\mathrm{cont}}
\newcommand{\sgn}{\mathrm{sgn}}
\newcommand{\Bock}{\mathrm{Bock}}
\newcommand{\tE}{\widetilde{E}}

\newcommand{\hT}{\widehat{T}}

\newcommand\blfootnote[1]{%
  \begingroup
  \renewcommand\thefootnote{}\footnote{#1}%
  \addtocounter{footnote}{-1}%
  \endgroup
}

\DeclareFontEncoding{OT2}{}{} 
\DeclareTextFontCommand{\textcyr}{\fontencoding{OT2}\fontfamily{wncyr}\fontseries{m}\fontshape{n}\selectfont}

\newcommand{\Sha}{\textcyr{Sh}}

\usepackage[textheight=220mm,textwidth=150mm]{geometry}

\title{On spectral sequences for semiabelian varieties over non-closed fields}

\author{Alexander Petrov and Alexei Skorobogatov}

\date{\today}
\begin{document}
\maketitle

\tableofcontents

\begin{abstract}
\noindent 
We give a new, short proof of the formula for the first potentially non-zero differential 
of the Hochschild--Serre spectral sequence for semiabelian varieties over non-closed fields.
We show that this differential is non-zero for the Jacobian of a curve 
when the image of the torsor of theta-characteristics
under the Bockstein map is non-zero.
An explicit example is a curve of genus 2 whose Albanese torsor is not divisible by 2.
When the Albanese torsor is trivial, we show that the Hochschild--Serre spectral sequence
for the Jacobian degenerates at the second page. We give a formula
for the differential of the Hochschild--Serre spectral sequence for a torus 
which computes its Brauer group.
Finally, we describe the differentials of the Hochschild--Serre spectral sequence for
a smooth projective curve, generalising a lemma of Suslin.
\end{abstract}

\blfootnote{Mathematics Subject Classification (2020): 18G40, 11G10}
\blfootnote{Keywords: spectral sequences, abelian varieties, torsors, tori, algebraic curves}

\section*{Introduction}

For an algebraic variety $X$ over a field $k$ with separable closure $k_\s$
and a prime $\ell$ not divisible by ${\rm char}(k)$, the action of the absolute Galois group $\Ga=\Gal(k_{\s}/k)$ on the \'etale cohomology groups $\H^i_{\et}(X_{k_{\s}},\Z_{\ell})$ contains arithmetic information about the variety $X$, such as the number of rational points in the case 
when $k$ is a finite field.

The action of the Galois group $\Ga$ on individual \'etale cohomology groups comes from its action on the \'etale cohomology complex $\RGamma_{\et}(X_{k_{\s}},\Z_{\ell})$. The phenomenon that we study in this paper is that the quasi-isomorphism class of this complex in general contains more information than individual cohomology Galois modules.

Remarkably, for a smooth projective variety $X$ over $k$, by Deligne's decomposition theorem \cite[Proposition 2.4]{deligne}, the rational \'etale cohomology complex $\RGamma_{\et}(X_{k_{\s}},\Q_{\ell})$ is quasi-isomorphic to the direct sum $\bigoplus\limits_{i\geq 0}\H^i_{\et}(X_{k_{\s}},\Q_{\ell})[-i]$ of individual cohomology groups as an object of the derived category of $\Q_{\ell}$-vector spaces equipped with an 
action of $\Gamma$.

This is no longer true for cohomology with integral coefficients, and is witnessed by the fact that the Hochschild--Serre spectral sequence 
\begin{equation}
E_{2}^{p,q}= \H^p(k,\H_\et^q(X_{k_\s},\Z_\ell))\Rightarrow\H_\et^{p+q}(X,\Z_\ell)
\label{ss ell}
\end{equation}
in general has non-zero differentials. 
This is the case when $\ell=2$, $k=\R$, $X(\R)=\emptyset$, as follows from work of
Artin--Verdier \cite[Corollary 2.6 (b)]{AV64}, see also \cite[\S 2]{Cox79}.
If $X$ is a smooth projective real curve with no real points and $\ell=2$, then
a differential on the 2nd (respectively, 3rd) page is non-zero if the genus 
of $X$ is odd (respectively, even), see Examples \ref{real curves odd genus} and \ref{real curves even genus}.

Our first main result is that (\ref{ss ell}) may fail to degenerate even for abelian varieties. We denote by $D(\Gamma,\Z/n)$ the derived category of the abelian category of discrete $\Z/n$-modules with a continuous action of $\Gamma=\Gal(k_{\s}/k)$.

\begin{inttheo}[{Theorem \ref{nonsplit abelian example}}]\label{intro: main example}
There exists a principally polarised abelian surface $A$ over $\Q$ such that $\RGamma_{\et}(A_{\qbar},\Z/2)$ is not quasi-isomorphic to 
$\bigoplus\limits_{i\geq 0}\H^i_{\et}(A_{\qbar},\Z/2)[-i]$ in the derived category 
$D(\Gamma,\Z/2)$, where 
$\Gamma=\Gal(\qbar/\Q)$. Specifically, the differential 
$$
\delta_2^{0,2}\colon\H^0(\Q,\H^2_{\et}(A_{\qbar},\Z/2))\to \H^2(\Q,\H^1_{\et}(A_{\qbar},\Z/2))$$
on the second page of the Hochschild--Serre spectral sequence for $A$ is non-zero.
\end{inttheo}

By Theorem \ref{jac} below, the Hochschild--Serre spectral sequence degenerates 
for any abelian variety
isomorphic to a direct factor of a product of Jacobians of curves each with a rational
divisor class of degree 1. In particular, the abelian surface $A$ from Theorem \ref{intro: main example} is an example of an abelian variety not of this type.

Let us sketch the proof of Theorem \ref{intro: main example}. Let $C$ be a smooth proper curve of genus $g$ over $k$
and let $J$ be the Jacobian of $C$. We identify $J$ with the dual abelian variety
using the canonical principal polarisation.
Using a general description of $\delta_2^{0,2}$ given in
Theorem \ref{intro: main semiabelian} below, we show that the
mod 2 first Chern class of the canonical principal polarisation is
an element of $\H^2_\et(J_{k_\s},\Z/2)^\Ga$ which
is sent by $\delta_2^{0,2}$ to the image of the class 
of the torsor of theta-characteristics under the Bockstein map 
$\H^1(k,J[2])\to \H^2(k,J[2])$. Thus $\delta_2^{0,2}\neq 0$ whenever
$[{\bf Pic}^{g-1}_{C/k}]\in \H^1(k,J)$ is not divisible by 2. 
A systematic method to construct such curves over number fields was found by B.~Creutz \cite{Cre13}.
The idea is to exploit the fine arithmetic structure available in this case,
namely, that $[{\bf Pic}^{g-1}_{C/k}]$ is automatically a 2-torsion element 
of the Tate--Shafarevich group $\Sha(J)$. 
The class $[{\bf Pic}^{g-1}_{C/k}]$ is divisible by 2 in $\H^1(k,J)$ if and only if its image
in $\Sha^2(k,J[2])$ under the Bockstein map is zero, which is equivalent to being
orthogonal to $\Sha^1(k,J[2])$ under the non-degenerate Poitou--Tate duality pairing. 
By compatibility of pairings, this happens if and only if
$[{\bf Pic}^{g-1}_{C/k}]$ is orthogonal to the image of $\Sha^1(k,J[2])$ in $\Sha(J)$
under the Cassels--Tate pairing.

We have two strategies to construct an element of $\Sha^1(k,J[2])$ whose image in $\Sha(J)$
pairs non-trivially with $[{\bf Pic}^{g-1}_{C/k}]$.

The first strategy, following Manin, interprets the Cassels--Tate pairing 
in terms of the Brauer--Manin obstruction on the everywhere locally soluble variety
$X={\bf Pic}^{g-1}_{C/k}$: any $t\in\Sha(J)$ gives rise to a Brauer class in
$\Br(X)$ whose Brauer--Manin pairing
with an arbitrary adelic point of $X$ equals the Cassels--Tate pairing of $t$ and $[X]$.
When $g=2$, the variety $X$ is the Albanese torsor of $C$, so $C$ is a closed subvariety of $X$.
If $C$ itself is everywhere locally soluble, we can take an adelic point of $X$ contained in
$C$ and then we only need to know the restriction of our Brauer class to $C$. 
The easiest way to arrange for this is to consider hyperelliptic curves $C$ given by
the equation $y^2=f(x)$, where the polynomial $f(x)\in k[x]$ of degree 6 is itself
everywhere locally soluble. 
In this setting, one can construct elements of $\Br(C)$ coming from 
$\Sha^1(k,J[2])\subset\H^1(k,J[2])$ by an explicit formula, 
which allows one to compute the Brauer--Manin pairing with the adelic points of $f(x)=0$. 
For a worked-out example 
borrowed from Creutz \cite[p.~941]{Cre13} and Creutz--Viray \cite[Theorem~6.7]{CV15}
see Theorem \ref{nonsplit abelian example} (a).

The second strategy is based on the approach of Poonen and Stoll \cite{PS}. 
The class $[{\bf Pic}^{g-1}_{C/k}]$
comes from the torsor of theta-characteristics in $\H^1(k, J[2])$, so if this torsor is
everywhere locally trivial and the curve $C$ is odd, which means that the Cassels--Tate
pairing of $[{\bf Pic}^{g-1}_{C/k}]$ with itself is non-zero, then 
we are done. For a worked-out example see Theorem \ref{nonsplit abelian example} (b).

Let us now give a general description of the first potentially non-trivial extension in the \'etale cohomology complex of a semiabelian variety over an arbitrary field $k$ (${\rm char}(k)\neq 2$).
It is easy to see that all differentials on the $i$-th page of (\ref{ss ell})
are zero for $\ell>i$, see Corollary \ref{1.3}. For the study of the differentials
on the second page we can thus assume that $\ell=2$.
For a $\Z_2$-module $M$ we denote by $Q(M)$ the module of quadratic functions on 
$M^{\vee}=\Hom_{\Z_2}(M,\Z_2)$ with values in $\Z_2$. It fits into the exact sequence
\begin{equation}\label{intro: quadratic ext2}
0\to M\to Q(M)\to M^{\otimes 2}\to \wedge^2 M\to 0,
\end{equation}
where the middle map sends a quadratic function $f:M^{\vee}\to\Z_2$ to the bilinear form
$$f(x+y)-f(x)-f(y)\in\Hom(M^{\vee}\otimes M^{\vee},\Z_2)\cong M^{\otimes 2}.$$ 
The sequence (\ref{intro: quadratic ext2}) is natural in $M$, so if $M$ is equipped with an action of $\Gamma$ then (\ref{intro: quadratic ext2}) becomes an exact sequence of $\Gamma$-modules. 

\begin{inttheo}[{Theorem \ref{p1}}]\label{intro: main semiabelian}
Let $A$ be a semiabelian variety over a field $k$
of characteristic not equal to $2$. The class in $\Ext^2_{\Ga}(\H^2_{\et}(A_{k_{\s}},\Z_2), \H^1_{\et}(A_{k_{\s}},\Z_2))$ corresponding to $\tau^{[1,2]}\RGamma_{\et}(A_{k_{\s}},\Z_2)$ is equal to the class of the extension 
$(\ref{intro: quadratic ext2})$ for $M=\H^1_{\et}(A_{k_{\s}},\Z_2)$. 
\end{inttheo}

In the case of abelian varieties this was proved in \cite{P}, but here we give a new shorter proof exhibiting an explicit quasi-isomorphism from the 2-term complex $Q(M)\to M^{\otimes 2}$ to the 
truncated bar complex of the Tate module $T_2(A)$. Theorem~\ref{intro: main semiabelian}, in particular, gives a description of the differentials 
$$\H^i(k, \H^2_{\et}(A_{k_{\s}},\Z_2))\to \H^{i+2}(k, \H^1_{\et}(A_{k_{\s}},\Z_2))$$
of (\ref{ss ell}), as well as the analogous differentials in the Hochschild--Serre spectral sequence with $\Z/2^n$-coefficients, for any $n$, see also Theorem \ref{p2}. 

Theorem \ref{intro: main example} shows that the class described in Theorem \ref{intro: main semiabelian} is not always zero for abelian varieties, but at the moment we do not know if it is ever non-zero for tori, cf. Remark \ref{tori: nonvanishing remark}.

In contrast, for the Jacobians of curves with a rational divisor class of degree $1$
the \'etale cohomology complex decomposes in all degrees. We give a proof of this result, which may
be well-known to the experts, but does not seem to be available in the literature.

\begin{inttheo}[{Theorem \ref{jacobian-decomposition}}] \label{jac}
Let $C_1,\ldots, C_m$ be geometrically
connected smooth proper curves over $k$ each admitting a $k$-rational divisor class of degree $1$,
for $i=1,\ldots,m$. Let $A$ be a direct factor of the product $\prod_{i=1}^m \Jac(C_i)$ of Jacobians of $C_1,\ldots,C_m$.
Let $n\geq 1$ be an integer not divisible by ${\rm char}(k)$.
Then $\RGamma_{\et}(A_{k_{\s}},\Z/n)$ is quasi-isomorphic to $\bigoplus\limits_{i\geq 0}\H^i_{\et}(A_{k_{\s}},\Z/n)[-i]$ in $D(\Gamma,\Z/n)$. In particular, the Hochschild--Serre spectral sequence $$E_2^{i,j}=\H^i(k,\H^j_{\et}(A_{k_{\s}},\Z/n))\Rightarrow \H^{i+j}_{\et}(A,\Z/n)$$ degenerates at the second page.
\end{inttheo}

Another class of examples of non-zero differentials in the Hochschild--Serre spectral sequence comes from varieties without a rational point. If a geometrically connected $k$-variety $X$ has a $k$-point, then $\RGamma_{\et}(X_{k_{\s}},\Z_{\ell})$ is quasi-isomorphic to $\Z_{\ell}\oplus\tau^{\geq 1}\RGamma_{\et}(X_{k_{\s}},\Z_{\ell})$ in the derived category of $\Gamma$-modules.
This is no longer true for varieties without a rational point, and we calculate the exact obstruction to splitting off $\H^0(X_{k_{\s}},\Z_{\ell})$ from $\tau^{\leq 1}\RGamma_{\et}(X_{k_{\s}},\Z_{\ell})$ when $X$ is a torsor for a semiabelian variety.
We denote by $\mathcal{A}(\Ga,\Z/n)$ 
the abelian category of discrete $\Z/n$-modules equipped with a continuous action of $\Ga$.

\begin{inttheo}[{Theorem \ref{1}}]\label{intro: torsors}
Let $X$ be a torsor for a semiabelian variety $A$ over a field $k$.
Let $n\geq 1$ be an integer not divisible by ${\rm char}\, k$.
Then the class in 
$$\Ext^2_{\mathcal{A}(\Ga,\Z/n)}(\H^1(X_{k_{\s}},\Z/n),\Z/n)\cong \Ext^2_{\mathcal{A}(\Ga,\Z/n)}(\H^1(A_{k_{\s}},\Z/n),\Z/n)\cong \H^2(k, A[n])$$
corresponding to $\tau^{[0,1]}\RGamma_{\et}(X_{k_{\s}},\Z/n)$ is equal to the image of the class of $X$ in $\H^1(k, A)$ under the Bockstein homomorphism $\H^1(k, A)\to \H^2(k, A[n])$.
\end{inttheo}

As we point out in Remark \ref{section remark}, in the setting of Theorem \ref{intro: torsors}, the obstruction to the existence of a rational $0$-cycle of degree $1$ coming from $\tau^{[0,1]}\RGamma_{\et}(X_{k_{\s}},\widehat{\Z})$ being non-split coincides with the obstruction coming from the failure  of the abelianised fundamental exact sequence for the arithmetic fundamental group to have a section.

We give several complements of the above results, summarised below. 

In Theorem \ref{tori: main gm coeffs} we describe a differential on the $3$rd page of the Hochschild--Serre spectral sequence with $\G_m$-coefficients for a {\em torus}, which allows one to 
compute the Brauer group of any algebraic torus, answering a question raised in \cite[p.~220]{CTS21}.
When the first version of this paper was completed, Julian Demeio told us about his result
that for a torus $T$ over a field $k$ of characteristic zero the natural map 
$\Br(T)\to\Br(T_{k_\s})^\Ga$ is surjective when $T$ is quasi-trivial or when
$k$ is a number field \cite[Theorem 1.1]{Dem}. This also follows from 
our Theorem \ref{tori: main gm coeffs},
see Corollary \ref{demeio}. Demeio's proof is based on an elaborate analysis of differentials
of the Hochschild--Serre spectral sequence by L.S.~Charlap, A.T.~Vasques, and C.H.~Sah (see the 
references in \cite{Dem}). The proof of Theorem \ref{tori: main gm coeffs} in this paper is self-contained: it reduces the computation of the differential with $\G_m$-coefficients to the case of finite coefficients which is handled in our explicit proof of Theorem \ref{intro: main semiabelian}. 

In Proposition \ref{curves: main} we combine Theorem \ref{intro: torsors} with Poincar\'e duality to describe all differentials on the second page of the Hochschild--Serre spectral sequence with coefficients in $\mu_n$ for a smooth proper {\em curve}. When all differentials on the second page vanish, in Proposition \ref{curves: 3rd page main} we describe the differentials on the $3$rd page of this spectral sequence, generalising a result of Suslin in the case of genus $0$.

\medskip

{\bf Notation.} For an object $C$ of the derived category $D(\mathcal{A})$ of an abelian category $\mathcal{A}$ its truncation $\tau^{[i,i+1]}C$ in any two consecutive degrees fits into a distinguished triangle
$$
\H^i(C)[-i]\to \tau^{[i,i+1]}C\to \H^{i+1}(C)[-i-1]\xrightarrow{\delta_i}\H^i(C)[-i+1]
$$
We refer to the map $$\delta_i\in \Hom_{D(\mathcal{A})}(\H^{i+1}(C)[-i-1],\H^i(C)[-i+1])\cong \Ext^2_{\mathcal{A}}(\H^{i+1}(C), \H^i(C))$$ as the extension class of $\tau^{[i,i+1]}C$.

For a profinite group $G$ and an integer $n$ we denote by $\mathcal{A}(G,\Z/n)$ 
the abelian category of discrete $\Z/n$-modules equipped with a continuous action of $G$, cf.~\cite[6.11]{Weibel}. We write $D(G,\Z/n)$ for the derived category of $\mathcal{A}(G,\Z/n)$.
For a prime $\ell$ we denote by $D(G,\Z_{\ell})$ the derived category of sheaves of $\Z_{\ell}$-modules on the pro-\'etale site $BG_{\proet}$, where $\Z_{\ell}$ is the sheaf of rings obtained from the topological ring $\Z_{\ell}$ with a trivial action of $G$, cf.~\cite[4.3]{bs}. Recall that an $\ell$-adically complete $\Z_{\ell}$-module $M$ equipped with a continuous action of $G$ gives rise to an object of the triangulated category $D(G,\Z_{\ell})$, and $\Hom_{D(G,\Z_{\ell})}(\Z_{\ell},M[i])$ is isomorphic to $i$-th continuous cohomology group $\H^i_{\cont}(G,M)$ defined using the complex of continuous cochains on $G$.

When $\Gamma=\Gal(k_{\s}/k)$ is the absolute Galois group of a field $k$, we denote $\H^i_{\cont}(\Gamma, M)$ by $\H^i(k, M)$ for any continuous $\Gamma$-module $M$.

\medskip

{\bf Acknowledgements.}
The first named author was supported by the Clay Research Fellowship and was a member of the Institute for Advanced Study in Princeton during the preparation of this paper.
Work on this paper benefited from visits of the second named author to the 
Institute for Advanced Study and the Max Planck Institute for Mathematics in Bonn.
We thank Jean-Louis Colliot-Th\'el\`ene and Tam\'as Szamuely
for help with references, and Vadim Vologodsky for explaining 
to us the result obtained in Theorem \ref{jac}.

\section{Semiabelian varieties}

Let $A$ be a semiabelian variety over $k$, that is, an extension of an abelian variety by a torus.
Let $\ell$ be a prime number, $\ell\neq{\rm char}(k)$. For all $i\geq 0$, $n\geq 1$, cup product induces an isomorphism 
$\H^i_\et(A_{k_\s},\Z/\ell^n)\simeq\wedge^i\H^1_\et(A_{k_\s},\Z/\ell^n)$
compatible with the natural action of the Galois group $\Ga=\Gal(k_{\s}/k)$. (By
\cite[Corollary VI.4.3]{EC} this follows from the same statement with $k_\s$
replaced by an algebraically closed field extension, which is proved in 
\cite[Lemma 4.1]{brionszamuely}.)

We would like to study the differentials 
$\delta_i^{p,q}\colon E_{i}^{p,q}\to E_i^{p+i,q-i+1}$
of the Hochschild--Serre spectral sequence
\begin{equation} \label{Leray1}
E_{2}^{p,q}= \H^p(k,\H_\et^q(A_{k_\s},\Z/\ell^n))\Rightarrow\H_\et^{p+q}(A,\Z/\ell^n).
\end{equation}

Let us make some simple remarks about these differentials.

\brem \label{1.1}
{\rm The origin of the group law on $A$ gives a section of the structure morphism
$\pi\colon A\to\Spec(k)$, hence the natural maps $\H^r_\et(k,\Z/\ell^n)\to\H^r_\et(A,\Z/\ell^n)$
are injective for all $r\geq 0$. Therefore, we have $\delta_i^{p,i-1}=0$ for all $i\geq 2$ and all $p\geq 0$.}
\erem

\ble \label{1.2}
Let $i\geq 2$, $p\geq 0$, $q\geq i-1$ be integers, and let $\ell$ be a prime. When $\ell-1$
divides $i-1$ we define $n=\val_\ell\big((i-1)/(\ell-1)\big)$.
\hfill $\Box$

\smallskip

{\rm (a)} Suppose $\ell\neq 2$. If $\ell-1$ does not divide $i-1$, then 
$\delta^{p,q}_i=0$.
Otherwise, we have $\ell^{\min\{q-i+1,n+1\}}\delta^{p,q}_i=0$.

{\rm (b)} Suppose $\ell=2$. If $i$ is even, then $2^{\min\{q-i+1,1\}}\delta^{p,q}_i=0$. If $i$ is odd,
then $2^{\min\{q-i+1,n+2\}}\delta^{p,q}_i=0$.
\ele
{\em Proof.} 
The spectral sequence (\ref{Leray1}) is functorial in $A$, so it is
compatible with multiplication by $m$ map $[m]\colon A\to A$,
for any integer $m$. The induced map on
$\H^r_\et(A_{k_\s},\Z/\ell^n)\cong\wedge^r\H^1_\et(A_{k_\s},\Z/\ell^n)$
is multiplication by $m^r$. Each $E_i^{p,q}$ for $i\geq 2$ is a subquotient of $E_{2}^{p,q}$, hence $[m]^*$ acts on it via multiplication by $m^q$. Since $[m]^*$ commutes with all differentials on the $i$-th page, we have $m^{q-i+1}\delta_{i}^{p,q}=\delta_{i}^{p,q}m^q$ for all $p,q\geq 0$. Thus we have 
$$(m^q-m^{q-i+1})\delta^{p,q}_i=m^{q-i+1}(m^{i-1}-1)\delta^{p,q}_i=0.$$
Taking $m=\ell$ we 
see that $\ell^{q-i+1}$ annihilates $\delta^{p,q}_i$. From now on we assume that
$(m,\ell)=1$. Then we have $(m^{i-1}-1)\delta^{p,q}_i=0$.

Suppose $\ell\neq 2$. If $\ell-1$ does not divide $i-1$, then we can find an integer $m$
such that $(m,\ell)=1$ and $m^{i-1}-1$ is non-zero
modulo $\ell$. Then $\delta^{p,q}_i=0$.

If $\ell-1$ divides $i-1$, we claim that the lowest value of $\val_\ell(m^{i-1}-1)$, where 
$(m,\ell)=1$, is $n+1$. If $(m,\ell)=1$, then
$m^{r\ell^n(\ell-1)}-1$ is clearly divisible by $\ell^{n+1}$, so it is enough to check that
$\val_\ell((1+\ell)^{r\ell^n(\ell-1)}-1)=n+1$ when $(r,\ell)=1$.
This is immediate for $n=0$, and the general case follows by induction in $n$.
This proves (a).

Suppose $\ell=2$. If $i$ is even, then taking $m=-1$ we prove the first statement of (b).
If $i$ is odd, so that $n\geq 1$, we claim that the smallest value of $\val_2(m^{r 2^n}-1)$, 
where $r$ and $m$ are odd, is $n+2$. We have
$$m^{r2^n}-1=(m^{2^n}-1)(m^{(r-1)2^n}+\ldots+1).$$
Since $r$ and $m$ are odd, the second factor in the right hand side is odd. Thus
we can assume that $r=1$. For $n=1$ the statement is obvious, and the general case
immediately follows by induction.
This proves (b). \hfill $\Box$

\medskip

Note that the statement of Remark \ref{1.1} is a particular case of Lemma \ref{1.2}.

\bco \label{1.3}
All differentials on the $i$-page of the spectral sequence $(\ref{Leray1})$ are zero for $\ell>i$.
\eco

Thus the differentials on the second page can be non-zero only for $\ell=2$, and
the first non-trivial case is the differential $\delta^{p,2}_2$, where $p\geq 0$.
Let us start with the spectral sequence
\begin{equation} \label{Leray2}
E_{2}^{p,q}= \H^p(k,\H_\et^q(A_{k_\s},\Z_2))\Rightarrow\H_\et^{p+q}(A,\Z_2).
\end{equation} 

Let $M$ be a free, finitely generated $\Z_2$-module. Recall that $S^2(M)$ is the quotient
of $M^{\otimes 2}$ by the $\Z_2$-submodule generated by $x\otimes y-y\otimes x$ for $x,y\in M$,
and $\wedge^2(M)$ is the quotient
of $M^{\otimes 2}$ by the $\Z_2$-submodule generated by $x\otimes x$ for $x\in M$.
We write $M^\vee=\Hom_{\Z_2}(M,\Z_2)$.
\bde
A function $f:M^{\vee}\to \Z_2$ is called quadratic if $f(x+y)-f(x)-f(y)$ is bilinear in $x$ and $y$.
Write $Q(M)$ for the $\Z_2$-module of quadratic functions on $M^{\vee}$. 
\ede

For example, $Q(\Z_2)$ is a free rank $2$ module consisting of functions of the form $a\cdot x+b\cdot \frac{x^2-x}{2}$ with $a,b\in\Z_2$.

There is a natural injection $M\to Q(M)$ sending
an element of $M$ to the linear function on $M^\vee$ that it defines.
The cokernel of this map is $\Hom_{\Z_2}(S^2(M^\vee),\Z_2)$.
One immediately checks (for example, by choosing a basis of $M$ and the dual basis of $M^\vee$)
that under the natural pairing 
$$M^{\otimes 2}\times (M^\vee)^{\otimes 2}\to\Z_2$$
the $\Z_2$-submodules 
$(M\otimes M)^{S_2}=\langle m\otimes m|m\in M\rangle$ and 
$\langle a\otimes b-b\otimes a|a,b\in M^{\vee}\rangle$
are exact annihilators of each other. Thus there is a canonical isomorphism
$\Hom_{\Z_2}(S^2(M^\vee),\Z_2)\cong(M\otimes M)^{S_2}$.
We obtain a canonical exact sequence 
\begin{equation}
0\to M\to Q(M)\to M^{\otimes 2}\to \wedge^2 M\to 0. \label{genius}
\end{equation}
The $\Z_2$-module $Q(M)$ contains the submodule of quadratic forms $M^\vee\to\Z_2$.
The map $Q(M)\to M^{\otimes 2}$ sends a quadratic form to the associated bilinear form.

The extension defined by (\ref{genius}) is equivalent to an extension with smaller terms. Denote by $(M^{\otimes 2})_{S_2,\sgn}=(M^{\otimes 2})/\langle m_1\otimes m_2+m_2\otimes m_1|m_1,m_2\in M\rangle$ the coinvariants of the involution $m_1\otimes m_2\mapsto -m_2\otimes m_1$ acting on $M^{\otimes 2}$. This module fits into the following exact sequence:
\begin{equation}
0\to M\stackrel{2}\lra M\to (M^{\otimes 2})_{S_2,{\rm sgn}}\to \wedge^2 M\to 0, \label{genius1}
\end{equation}
where the map $M\to (M^{\otimes 2})_{S_2,\sgn}$ sends $m$ to $m\otimes m$, and the rightmost map sends $m_1\otimes m_2\in (M^{\otimes 2})_{S_2,\sgn}$ to $m_1\wedge m_2\in \wedge^2 M$. A direct computation shows:
\ble 
There is a commutative diagram given by
$$
\begin{tikzcd}
M\arrow[r]\arrow[d, equal] & Q(M)\arrow[r]\arrow[d] & M^{\otimes 2}\arrow[r]\arrow[d, two heads] & \wedge^2 M\arrow[d, equal] \\
M\arrow[r, "2"] & M\arrow[r] & (M^{\otimes 2})_{S_2,\sgn}\arrow[r] & \wedge^2 M
\end{tikzcd}
$$
where the map $Q(M)\to M$ sends a quadratic function $f$ on $M^{\vee}$ to the linear function $m\mapsto 4f(m)-f(2m)$.
\ele
 
The extension (\ref{genius1}) is the Yoneda product of 
\begin{equation}
0\to M\stackrel{2}\lra M\to M/2\to 0 \label{Bock}
\end{equation}
and 
\begin{equation}
0\to M/2\to (M^{\otimes 2})_{S_2,{\rm sgn}}\to \wedge^2 M\to 0. \label{alpha*}
\end{equation}
All of the above constructions and exact sequences are compatible with the action of 
the automorphisms of the $\Z_2$-module $M$. In particular, if $M$ is a $\Z_2$-module equipped with the $2$-adic topology and a continuous action of $\Gamma$, all the sequences above are naturally sequences of $\Gamma$-modules.

\medskip

The following statement is the particular case $p=2$ of \cite[Corollary 9.5 (1)]{P}. Here we give a short elementary proof of this result.

\bthe \label{p1}
Let $k$ be a field of characteristic different from $2$, and let $A$ be a semiabelian variety over $k$. 
Write $M=\H^1_\et(A_{k_\s},\Z_2)$ so that
$\wedge^i M\cong\H^i_\et(A_{k_\s},\Z_2)$ are continuous $\Ga$-modules for $i\geq 0$.
The differential 
$$\delta_2^{p,2}\colon \H^p(k,\wedge^2 M)\to \H^{p+2}(k,M)$$
of the spectral sequence $(\ref{Leray2})$
equals the connecting maps of the equivalent $2$-extensions of $\Ga$-modules $(\ref{genius})$
and $(\ref{genius1})$.
For any $n\geq 1$, the differential 
$$\delta_2^{p,2}\colon \H^p(k,\wedge^2(M/2^n))\to \H^{p+2}(k,M/2^n)$$
of the spectral sequence $(\ref{Leray1})$ with $\ell=2$
equals the connecting map of the $2$-extension of $\Ga$-modules
\begin{equation}
0\to M/2^n\to Q(M)/2^n\to (M/2^n)^{\otimes 2}\to \wedge^2 (M/2^n)\to 0. \label{mod2^n}
\end{equation}
\ethe
{\em Proof.}
We begin by recasting \'etale cohomology of $A_{k_{\s}}$ in terms of group cohomology of the profinite abelian group $M^{\vee}$, using that semiabelian varieties are \'etale $K(\pi,1)$-spaces:
\ble \label{odin}
The \'etale cohomology complex $\RGamma_{\et}(A_{k_{\s}},\Z_2)$ is identified with the continuous group cohomology complex $\RGamma_{\cont}(M^{\vee},\Z_2)$ as an object of $D(\Gamma,\Z_2)$, where the action of $\Gamma$ on the latter complex is induced by its action on the abelian group $M$. Likewise, the \'etale cohomology complex $\RGamma_{\et}(A_{k_{\s}},\Z/2^n)\in D(\Gamma,\Z/2^n)$ with finite coefficients is identified with $\RGamma_{\cont}(M^{\vee},\Z/2^n)\simeq  \RGamma_{\cont}(M^{\vee},\Z_2)\otimes^L_{\Z_2}\Z/2^n$.
\ele

\brem{\rm \label{rmk}
It follows from Lemma \ref{odin} that canonical isomorphisms of $\Ga$-modules
$\H^i_{\rm cont}(M^\vee,\Z_2)\cong \wedge^i M$ for all $i\geq 0$ identify
the spectral sequence $(\ref{Leray2})$ with the Hochschild--Serre 
spectral sequence for continuous group cohomology of the semi-direct product $M^\vee\rtimes\Ga$ 
for the natural action of $\Ga$ on the Tate module $T_2(A_{k_{\s}})\cong M^\vee$, with values in 
the trivial module $\Z_2$:
\begin{equation}
E_2^{p,q}=\H^p_{\rm cont}(\Ga,\H^q_{\rm cont}(M^\vee,\Z_2))\Rightarrow \H^{p+q}_{\rm cont}(M^\vee\rtimes\Ga,\Z_2).
\label{contHS}
\end{equation}
}
\erem
Given Lemma \ref{odin}, Theorem \ref{p1} follows from the following description of the extension class in $\tau^{[1,2]}\RGamma(M^{\vee},\Z_2)$:

\bpr \label{dva}
The extension class in $\Ext^2_{\Gamma}(\wedge^2 M, M)$ defined by $\tau^{[1,2]}\RGamma(M^{\vee},\Z_2)\in D(\Gamma,\Z_2)$ can be represented by (\ref{genius}). 
The extension class in $\Ext^2_{\Gamma}(\wedge^2(M/2^n),M/2^n)$ defined by $\tau^{[1,2]}\RGamma(M^{\vee},\Z/2^n)\in D(\Gamma,\Z/2^n)$ is represented by the mod $2^n$ reduction (\ref{mod2^n}) of  (\ref{genius}).
\epr

{\em Proof of Lemma \ref{odin}.} 
For any connected scheme $X$ with a geometric base point $\bar{x}$ and any integer $m\geq 1$
there is a natural map
\begin{equation}\label{pi1coh to etalecoh map}
\RGamma(\pi_1^{\et}(X,\bar{x}),\Z/m)\to \RGamma_{\et}(X,\Z/m)
\end{equation} 
induced by the pullback along the map of sites $X_{\et}\to X_{\mathrm{f\acute et}}$, where the target site consists of schemes finite \'etale over $X$, cf.~\cite[\S 5]{olsson}. Moreover, if $X=Y_{k_{\s}}$, where $Y$ is a scheme over $k$, and $\bar{x}$ is above a $k$-point of $Y$, this map naturally lifts to a map in $D(\Gamma,\Z/m)$, where the action of $\Gamma$ on the continuous cochain complex $\RGamma(\pi_1^{\et}(X,\bar{x}),\Z/m)$ is induced by the Galois action on the fundamental group $\pi_1^{\et}(X,\bar{x})$.

The map (\ref{pi1coh to etalecoh map}) induces a natural map
\begin{equation}
\H^1_{\rm cont}(\pi_1^{\et}(X,\bar{x}),\Z/m)\to \H^1(X,\Z/m). \label{25a}
\end{equation}
Since \'etale $X$-torsors for $\Z/m$ are the same as Galois covers of $X$ with Galois group $\Z/m$,
this map is an isomorphism \cite[Expos\'e XI, \S 5]{SGA1}.
Consider the composition
\begin{equation}\label{pi1coh to etalecoh map prime to p}
\RGamma(\pi_1^{\et}(X,\bar{x})^{(p')},\Z/m)\to \RGamma(\pi_1^{\et}(X,\bar{x}),\Z/m)\to \RGamma_{\et}(X,\Z/m),
\end{equation} 
where the first map is induced by the surjection $\pi_1^{\et}(X,\bar{x})\twoheadrightarrow\pi_1^{\et}(X,\bar{x})^{(p')}$ onto the prime-to-$p$ completion of the fundamental group where $p={\rm char}(k)$; if ${\rm char}(k)=0$ 
we set $\pi_1^{\et}(X,\bar{x})^{(p')}:=\pi_1^{\et}(X,\bar{x})$.
The map (\ref{pi1coh to etalecoh map prime to p}) induces an isomorphism on $\H^1$ 
whenever $m$ is not divisible by ${\rm char}(k)$.

Let us now specialise to the case $X=A_{k_{\s}}$, where $A$ is a semiabelian variety over $k$.
We take $\bar x$ to be the origin of the group law on $A$ and omit it from the notation.

We know that the cohomology ring of the target of (\ref{pi1coh to etalecoh map prime to p}) 
is the free exterior algebra on $\H^1(A_{k_\s},\Z/m)$. 
We claim that the same is true for the source of 
(\ref{pi1coh to etalecoh map prime to p}). Indeed, \cite[Proposition 1.1 b)]{brionszamuely} implies 
that $\pi_1^{\et}(A_{k_{\s}})^{(p')}$ is a finite product of copies of $\widehat{\Z}^{(p')}$.
The cohomological dimension of $\widehat{\Z}^{(p')}$ is 1 and the cohomology groups
$\H^0(\widehat{\Z}^{(p')},\Z/m)\simeq \H^1(\widehat{\Z}^{(p')},\Z/m)\simeq\Z/m$ 
are free $\Z/m$-modules, hence the claim follows from the K\"unneth formula. Thus
for $X=A_{k_{\s}}$ the map (\ref{pi1coh to etalecoh map prime to p}) is a quasi-isomorphism 
in $D(\Gamma,\Z/m)$.

Next, $\pi_1^{\et}(A_{k_{\s}})^{(p')}$ is the product of the $2$-adic Tate module 
$T_2(A_{k_{\s}})\cong M^{\vee}$ and a 2-divisible profinite group.
Passing to the limit over $m=2^n$ we obtain a quasi-isomorphism 
$\RGamma(\pi_1^{\et}(A_{k_{\s}})^{(p')},\Z_2)\simeq \RGamma(M^{\vee},\Z_2)$
in $D(\Gamma,\Z_2)$, hence also a quasi-isomorphism 
$\RGamma_{\et}(A_{k_{\s}},\Z_2)\simeq \RGamma(M^{\vee},\Z_2)$.
\hfill$\Box$

\medskip

{\em Proof of Proposition \ref{dva}.}
We can calculate $\RGamma(M^{\vee},\Z_2)$ using the standard bar complex 
of $\Z_2$-modules for computing continuous cohomology of the group $M^{\vee}$:
$$
0\lra\Z_2\stackrel{d_0=0}\lra {\rm Func}(M^\vee,\Z_2)\stackrel{d_1}\lra 
{\rm Func}(M^\vee\times M^\vee,\Z_2)\stackrel{d_2}\lra\ldots,
$$
where ${\rm Func}$ denotes the $\Z_2$-module of continuous functions.
Considering this complex as a complex of $\Ga$-modules for the natural action of the Galois group 
$\Ga$ induced by its action on $M$, we get a complex in $D(\Gamma,\Z_2)$ that calculates
the cohomology of the semi-direct product $M^\vee\rtimes\Ga$.

The differential $d_1$ sends a function $f\colon M^\vee\to\Z_2$ to the function
$$(x,y)\mapsto f(x+y)-f(x)-f(y).$$ The inclusion $Q(M)\subset {\rm Func}(M^\vee,\Z_2)$
gives rise to a commutative diagram
\begin{equation}
\begin{split}
\label{main thm proof commdiag}\xymatrix{Q(M)\ar[r]\ar[d]&M^{\otimes 2}\ar[d]\\
{\rm Func}(M^\vee,\Z_2)\ar[r]^{d_1 \ \ \ \ }&{\rm Func}(M^\vee\times M^\vee,\Z_2)}
\end{split}
\end{equation}
The right-hand vertical map is the inclusion of bilinear functions on $M^\vee\times M^\vee$
into all continuous functions. Commutativity of the diagram is immediate from the definitions of the maps. The submodule $M\subset Q(M)$ of linear functions maps isomorphically to $\H^1_{\cont}(M^{\vee},\Z_2)$. Moreover, the right vertical map lands in the kernel of the next differential $d_2$, and the induced map $M^{\otimes 2}\to \H^2_{\cont}(M^{\vee},\Z_2)$ coincides with the cup-product map $\H^1_{\cont}(M^{\vee},\Z_2)^{\otimes 2}\to \H^2_{\cont}(M^{\vee},\Z_2)$.

Exactness of (\ref{genius}) then implies that (\ref{main thm proof commdiag}) gives a quasi-isomorphism in $D(\Gamma,\Z_2)$ between the two-term complex $Q(M)\to M^{\otimes 2}$
and the truncation $\tau^{[1,2]}\RGamma(M^{\vee},\Z_2)$, as desired. The description of the analogous truncation $\tau^{[1,2]}\RGamma(M^{\vee},\Z/2^n)$ now follows formally because cohomology groups $\H^i_{\cont}(M^{\vee},\Z_2)$ as well as terms of the extension (\ref{genius}) are torsion-free. \hfill$\Box$

\medskip

Theorem \ref{p1} is now proved because, by Lemma \ref{odin}, the differentials under consideration in spectral sequence (\ref{Leray2}) (resp. (\ref{Leray1})) are induced by the $\Ext^2$ class of $\tau^{[1,2]}\RGamma(M^{\vee},\Z_2)$ (resp. $\tau^{[1,2]}\RGamma(M^{\vee},\Z/2^n)$).  

\medskip

The rest of this section is mainly devoted to the proof of Theorem \ref{p2} below, which gives a more computable description of the differential $\delta_2^{p,2}$
of the spectral sequence $(\ref{Leray1})$ with coefficients $\Z/\ell^m$.

When $N$ is a free, finitely generated
$\Z_2$-module with continuous action of $\Ga$, we denote by Bock
the morphism $N/2^m\to N[1]$ defined by the exact sequence
$$0\to N\stackrel{2^m}\lra N\to N/2^m\to 0.$$
For $i\geq 1$, by a slight abuse of notation, we also denote by Bock the composition 
$N/2^m\to N[1]\to (N/2^i)[1]$. This map is
defined by the exact sequence 
$$0\to N/2^i\to N/2^{m+i}\to N/2^m\to 0,$$
which is the push-out of the previous exact sequence by the map $N\to N/2^i$.

Let us also introduce the following notation. For an $\F_2$-vector space $V$, equipped with an action of $\Gamma$ continuous with respect to the discrete topology on $V$, consider the short exact sequence
\begin{equation}\label{alpha}
0\to V\xrightarrow{v\mapsto v\cdot v} S^2 V\xrightarrow{v_1\cdot v_2\mapsto v_1\wedge v_2} \wedge^2 V\to 0.
\end{equation}
Its extension class is a map $\wedge^2 V\to V[1]$ in $D(\Gamma,\F_2)$ which we denote by $\alpha(V)$. Note that the extension (\ref{alpha*}) is obtained from (\ref{alpha}) for $V=M/2$ by pulling back along $\wedge^2 M\to \wedge^2(M/2)$. 
We have $S^2V=(V^{\otimes 2})_{S_2}\cong (V^{\otimes 2})_{S_2,{\rm sgn}}$. The following lemma shows that when the dimension of $V$ is small,
the connecting map $\alpha(V)$ defined by (\ref{alpha}) is zero.

\ble \label{small}
{\rm (i)} If the rank of $V$ is $2$, then $(\ref{alpha})$ is split.

{\rm (ii)}
If the rank of $V$ is $3$, then the connecting map $(\wedge^2V)^\Ga\to\H^1(\Ga,V)$
defined by $(\ref{alpha})$ is zero.
\ele
{\em Proof.} (i) Let $u,v,w$ be the three non-zero elements of $V$. 
Then the unique non-zero element
in $\wedge^2 V$ lifts to
$uv+vw+wu\in  S^2V$. This lifting is $\GL(2,\F_2)$-equivariant. 

(ii) Assume that there is a non-zero element $x\in (\wedge^2 V)^\Ga$. 
We have a perfect bilinear pairing of $\F_2$-vector spaces
$V\times \wedge^2V\to \F_2$. The elements of $V$ that pair trivially with $x$
form an $\F_2$-subspace $N\subset V$ of dimension 2
which is stable under the action of $\Ga$. Moreover, $x\in \wedge^2N$.
We have a commutative diagram of $\Ga$-modules with exact rows
$$\xymatrix{0\ar[r]&N\ar[r]\ar[d]&S^2 N\ar[r]\ar[d]&\wedge^2 N\ar[r]\ar[d]&0\\
0\ar[r]& V\ar[r]& S^2V\ar[r]& \wedge^2 V\ar[r]& 0}$$
By part (i), the top sequence is split. Thus $x$ lifts to an element of $(S^2V)^\Ga$.
\hfill $\Box$

\ble \label{permutational}
{\rm (i)} Suppose that $V$ is a permutational $\Ga$-module, that is, the action of $\Ga$ on $V$
preserves an $\F_2$-basis of $V$. Then $(\ref{alpha})$ is split.

{\rm (ii)} Suppose that $k=\R$. Then $(\ref{alpha})$ is split.
\ele
{\em Proof.} (i) Let $e_1,\ldots,e_n$ be a $\Ga$-stable basis of $V$. The union of sets 
$\{e_i^2\}$ and $\{e_ie_j|i<j\}$ is a $\Ga$-stable basis of $S^2V$, 
which gives a $\Ga$-equivariant splitting of $(\ref{alpha})$.

(ii) By (i) it is enough to show that any representation of $\Z/2$ in $V$ is permutational,
but this is well-known. Indeed,
endow $V$ with the structure of an $\F_2[x]$-module by letting $x$ 
act as the generator of $\Z/2$. By the classification of finitely generated torsion modules over a PID,
the $\F_2[x]$-module $V$ is a direct sum of cyclic submodules isomorphic to
$\F_2=\F_2[x]/(x-1)$ or $\F_2[x]/(x^2-1)$, both of which are permutational. \hfill $\Box$

\medskip

The proof of Theorem \ref{p2} relies on the following general computation.

\ble\label{reduction of bock composition}
Let $p$ be a prime. Let $X, Y$ be free, finitely generated $\Z_p$-modules equipped with a continuous action of $\Gamma$. Suppose we are given an extension 
$$0\to X/p\to E\to Y/p\to 0$$ of $\F_p$-vector spaces equipped with an action of $\Gamma$.
Let $\delta_E\colon Y/p\to X/p[1]$ be the corresponding map in $D(\Gamma,\F_p)$. Let $\tE$ be the $\Z_p$-module with $\Gamma$-action obtained by pulling back $E\twoheadrightarrow Y/p$ along $Y\twoheadrightarrow Y/p$. The exact sequence 
$$0\to X\xrightarrow{p}X\xrightarrow{\psi} \tE\to Y\to 0,$$ 
where $\psi$ is the composition $X\to X/p\hookrightarrow \tE$,
defines a map $\beta_E:Y\to X[2]$ in $D(\Gamma,\Z_p)$.
Then the reduction of $\beta_E$ modulo $p^m$ is equal to the difference of the composed maps (\ref{deltabock}) and (\ref{bockdelta}):
\begin{equation}\label{deltabock}
Y/p^m\xrightarrow{\Bock_{Y/p^{m+1}}} Y/p[1]\xrightarrow{\delta_E[1]} X/p[2]\xrightarrow{}X/p^m[2],
\end{equation}
\begin{equation}\label{bockdelta}
Y/p^m\to Y/p\xrightarrow{\delta_E} X/p[1]\xrightarrow{\Bock_{X/p^{m+1}}[1]}X/p^m[2].
\end{equation}
\ele
{\em Proof.} The mod $p^m$ reduction of the map $Y\to X[2]$ is the Yoneda class of the complex $(X\xrightarrow{\psi} \tE)\otimes^L_{\Z_p}\Z/p^m$ concentrated in degrees $-1$ and $0$. It can be computed as the totalisation of the bicomplex 
\begin{equation}\label{modp bicomplex}
\begin{tikzcd}
X\arrow[r,"{\psi}"] & \tE \\
X\arrow[r,"{\psi}"]\arrow[u, "{p^m}"] & \tE\arrow[u, "{p^m}"]
\end{tikzcd}
\end{equation}
The compositions of (\ref{deltabock}) and (\ref{bockdelta}) are represented, respectively, by the following Yoneda extensions
\begin{equation}\label{deltabock yoneda}
0\to X/p^m\xrightarrow{a} R_m\xrightarrow{b} Y/p^{m+1}\xrightarrow{c} Y/p^m\to 0
\end{equation}
\begin{equation}\label{bockdelta yoneda}
0\to X/p^m\xrightarrow{a'} X/p^{m+1}\xrightarrow{b'} T_m\xrightarrow{c'} Y/p^m\to 0,
\end{equation}
where $R_m$ (respectively, $T_m$) is the following pushout (respectively, pullback)
$$
\begin{tikzcd}
X/p\arrow[r]\arrow[d] & E\arrow[d] & & T_m \arrow[d]\arrow[r] & Y/p^m\arrow[d] \\
X/p^{m}\arrow[r] & R_m & & E \arrow[r] & Y/p
\end{tikzcd}
$$
and the map $b:R_m\to Y/p^{m+1}$ is the composition $R_m\to Y/p\hookrightarrow Y/p^{m+1}$, while $b':X/p^{m+1}\to T_m$ is the composition $X/p^{m+1}\twoheadrightarrow X/p\to T_m$. The result of subtracting the map $Y/p^m\to X/p^m[2]$ corresponding to (\ref{deltabock yoneda}) from that corresponding to (\ref{bockdelta yoneda}) is then represented by the extension
\begin{equation}\label{difference yoneda}
X/p^m\xrightarrow{(a,0)} \frac{R_m\oplus X/p^{m+1}}{(a(x),-a'(x))|x\in X/p^m}\xrightarrow{b\oplus b'} \ker(Y/p^{m+1}\oplus T_m\xrightarrow{c-c'} Y/p^m)\xrightarrow{c} Y/p^m    
\end{equation}
We can now write down a chain level map from the totalisation of (\ref{modp bicomplex}) to the two-term complex from (\ref{difference yoneda}):
\begin{equation}\label{bock main map}
\begin{tikzcd}
X\arrow[r,"{(\psi,-p^m)}"] & \tE\oplus X \arrow[r,"{(p^m,\psi)}"]\arrow[d,"{f_1}"] & \tE\arrow[d, "{f_2}"] \\
& \frac{R_m\oplus X/p^{m+1}}{(a(x),-a'(x))|x\in X/p^m}\arrow[r,"{b\oplus b'}"] & \ker(Y/p^{m+1}\oplus T_m\xrightarrow{c-c'} Y/p^m)
\end{tikzcd}
\end{equation}
Let $f_1$ be the map induced by the direct sum of maps $\tE\to E\to R_m$ and $X\twoheadrightarrow X/p^{m+1}$, 
and $f_2$ is the sum of compositions $\tE\to Y\to Y/p^{m+1}$ and $\tE\twoheadrightarrow T_m$. The map $f_2$ indeed lands in the kernel of $c-c'$, the square in (\ref{bock main map}) commutes, and the composition $f_1\circ (\psi,-p^m)$ is zero.

The complexes given by rows of (\ref{bock main map}) both have two non-zero cohomology groups given by $X/p^m$ and $Y/p^m$. The map of complexes induced by $f_1,f_2$ induces the identity map on cohomology, giving the desired identification between the two classes in 
$\Ext^2_{\mathcal{A}(\Ga,\Z_p)}(Y/p^m,X/p^m)$. 
\hfill $\Box$

\medskip

We now specialise the above computation to our case of interest. As before, for our semiabelian variety $A$ we denote the continuous $\Gamma$-module $\H^1_{\et}(A_{k_{\s}},\Z_2)$ by $M$. Recall that we write $\alpha(M/2)\colon \wedge^2 (M/2)\to (M/2)[1]$ for the morphism in $D(\Gamma,\F_2)$ (as well as the corresponding morphism between the same objects in $D(\Gamma,\Z_2)$) defined by (\ref{alpha}) with $V=M/2$.

The next statement is a more general variant of the particular case $p=2$ of \cite[Corollary~9.5 (2)]{P}.
We give an elementary proof of this result using Lemma \ref{reduction of bock composition}.

\bthe \label{p2}
Let $k$ be a field of characteristic different from $2$. 
Let $A$ be a semiabelian variety over $k$. Consider the spectral sequence
$$E_{2}^{p,q}= \H^p(k,\H_\et^q(A_{k_\s},\Z/2^m))\Rightarrow \H_\et^{p+q}(A,\Z/2^m).$$
The differential 
$\delta_2^{p,2}\colon \H^p(k,\wedge^2 (M/2^m))\to \H^{p+2}(k,M/2^m)$
is given by
$$
\alpha({M/2})\circ \Bock^{m,1}_{\wedge^2 M}-\Bock^{1,m}_M\circ\alpha(M/2),
$$
where $\Bock^{m,1}_{\wedge^2 M}:\H^p(k,\wedge^2(M/2^m))\to \H^{p+1}(k,\wedge^2(M/2))$ and $\Bock^{1,m}_M: \H^{p+1}(k, M/2)\to \H^{p+2}(k, M/2^m)$ are Bockstein homomorphisms.
\ethe
{\em Proof.} By Theorem \ref{p1} the degree $2$ extension class corresponding to $\tau^{[1,2]}\RGamma_{\et}(A_{k_\s},\Z_2)$ is equal to the composition
$$
\wedge^2M\to\wedge^2(M/2)\xrightarrow{\alpha(M/2)}M/2[1]\xrightarrow{\Bock_M}M[2]
$$
To calculate the induced map on mod $2^m$ reductions, we apply 
Lemma \ref{reduction of bock composition} to the case $p=2$ with $X=M$, $Y=\wedge^2 M$, and the extension $E$ given by 
$$0\to M/2\to S^2(M/2)\to \wedge^2(M/2)\to 0.$$ 
We obtain that the extension class of $\tau^{[1,2]}\RGamma_{\et}(A_{k_{\s}},\Z/2^m)$ is equal to the difference of the following compositions:
\begin{equation}
\wedge^2(M/2^m)\lra \wedge^2(M/2)\stackrel{\alpha(M/2)}\lra (M/2)[1]
\stackrel{\rm Bock}\lra(M/2^m)[2], \label{genius2}
\end{equation}
\begin{equation}
\wedge^2(M/2^m)\stackrel{\rm Bock}\lra \wedge^2(M/2)[1]\stackrel{\alpha(M/2)[1]}\lra (M/2)[2]
\lra(M/2^m)[2]. \label{genius3}
\end{equation}
The differential $\delta_2^{p,2}$ is equal to the map induced by this difference on the degree $p$ derived functor of $\Gamma$-invariants, which gives the desired result. \hfill$\Box$

\section{Jacobians}

Let $C$ be a smooth, projective, geometrically integral curve over a field $k$. 
Let ${\bf Pic}_{C/k}$ be the Picard scheme of $C$.
Denote by $J:={\bf Pic}^0_{C/k}$ the Jacobian of $C$.
The connected component  ${\bf Pic}_{C/k}^{1}\subset{\bf Pic}_{C/k}$ 
parametrising divisors of degree 1 is 
a torsor for $J$, usually called the Albanese torsor. There is a
canonical map $C\hookrightarrow {\bf Pic}_{C/k}^{1}$. After choosing
a point $x_0\in C(k_\s)$, this map is identified with the usual Abel--Jacobi map 
$C_{k_\s}\to J_{k_\s}$ sending $x$ to the class of the divisor $x-x_0$.

It is well known that
${\bf Pic}_{C/k}(k)$ is canonically isomorphic to $\Pic(C_{k_\s})^\Ga$,
see, e.g., \cite[Corollary 2.5.9]{CTS21}. Thus the class
of the torsor ${\bf Pic}^1_{C/k}(k)$ in $\H^1(k,J)$ is zero if and only if $C$ has
a $k$-rational divisor class of degree 1.

Let $n$ be an integer not divisible by ${\rm char}(k)$. As before, we denote by $D(\Gamma,\Z/n)$ the derived category of the abelian category of discrete $\Z/n$-modules equipped with a continuous action of $\Gamma=\Gal(k_\s/k)$.

Suppose that $A$ and $B$ are abelian varieties over $k$ such that there are quasi-isomorphisms
$$\RGamma_{\et}(A_{k_\s},\Z/n)\simeq \bigoplus\limits_{i\geq 0} \H^i_{\et}(A_{k_\s},\Z/n)[-i],
\quad 
\RGamma_{\et}(B_{k_\s},\Z/n)\simeq \bigoplus\limits_{j\geq 0} \H^j_{\et}(B_{k_\s},\Z/n)[-j]$$
in $D(\Gamma,\Z/n)$.
The K\"unneth formula gives a quasi-isomorphism
$$\RGamma_{\et}(A_{k_\s}\times B_{k_\s},\Z/n)\simeq
\RGamma_{\et}(A_{k_\s},\Z/n)\otimes^L_{\Z/n}\RGamma_{\et}(B_{k_\s},\Z/n).$$
The groups $\H^i_{\et}(A_{k_\s},\Z/n)$ and $\H^j_{\et}(B_{k_\s},\Z/n)$ are free $\Z/n$-modules
for all $i$ and $j$,
thus for each $m\geq 0$ we have isomorphisms
$$\H^m_\et(A_{k_\s}\times B_{k_\s})\cong\bigoplus_{i+j=m}
\H^i_{\et}(A_{k_\s},\Z/n)\otimes_{\Z/n}\H^j_{\et}(B_{k_\s},\Z/n).$$
We deduce a quasi-isomorphism
$$\RGamma_{\et}(A_{k_\s}\times B_{k_\s},\Z/n)\simeq \bigoplus\limits_{m\geq 0}
 \H^m_{\et}(A_{k_\s}\times B_{k_\s},\Z/n)[-m].$$
Conversely, given such a quasi-isomorphism, we can produce a quasi-isomorphism
\begin{multline}\RGamma_{\et}(A_{k_{\s}},\Z/n)\xrightarrow{p_1^*}\RGamma_{\et}(A_{k_\s}\times B_{k_\s},\Z/n)\simeq \bigoplus\limits_{m\geq 0}
 \H^m_{\et}(A_{k_\s}\times B_{k_\s},\Z/n)[-m]\\ \xrightarrow{(\id_{A_{k_{\s}}}\times e_B)^*}\bigoplus\limits_{m\geq 0}\H^m_\et(A_{k_{\s}},\Z/n)[-m].
\end{multline}
Hence direct products of abelian varieties with decomposable \'etale cohomology complex, as well as direct factors of such varieties, also have decomposable \'etale cohomology complexes. Therefore Theorem \ref{jac} of the introduction is a consequence of the following

\bthe\label{jacobian-decomposition} 
Let $C$ be a smooth, projective, geometrically integral curve over a field $k$. 
Let $n$ be an integer not divisible by ${\rm char}(k)$.
If $C$ has a $k$-rational divisor class of degree $1$, then there exists a quasi-isomorphism
\begin{equation}\label{jacobian-dec-formula}
\RGamma_{\et}(J_{k_\s},\Z/n)\simeq \bigoplus\limits_{i\geq 0} \H^i_{\et}(J_{k_\s},\Z/n)[-i]
\end{equation}
in $D(\Gamma,\Z/n)$. In particular, the spectral sequence
$$E^{p,q}_2=\H^p(k,\H^q_\et(J_{k_\s},\Z/n))\Rightarrow\H^{p+q}_\et(J,\Z/n)$$ 
degenerates at the second page, and for all $r\geq 0$ there are isomorphisms
\begin{equation}\label{jacobian-dec-oncoh-formula}
\H^r_{\et}(J,\Z/n)\cong\bigoplus\limits_{i+j=r}\H^i(k, \H^j_{\et}(J_{k_\s},\Z/n)).
\end{equation}
\ethe

{\em Proof.} We will deduce this decomposition from the fact that $J$ admits integral K\"unneth projectors defined over the field $k$, as proved in \cite[Theorem 1.4]{Suh}.

Composing the canonical map $C\to {\bf Pic}_{C/k}^{1}$ with an isomorphism of torsors
${\bf Pic}_{C/k}^{1}\cong J$ provided by the given degree $1$ divisor class, we get a map
$\alpha:C\to J$. Denote by $\alpha_n:\Sym^nC\to J$ the maps from the symmetric powers of $C$ induced by $\alpha$ using the group structure on $J$. For each $1\leq n\leq g$, denote by $w^{[n]}\in Z^n(J)$ the codimension $n$ cycle in the Jacobian variety obtained by pushing forward the fundamental cycle $[\Sym^{g-n} C]$ along $\alpha_{g-n}$. We also denote by $w^{[0]}\in Z^0(J)$ the fundamental cycle of $J$ itself.

Following \cite[Theorem 4.2.3]{Suh}, for each $0\leq i\leq 2g$ we define a codimension $g$ cycle on $J\times J$ by the formula
$$
\pi_i:=(-1)^i\displaystyle\sum_{\substack{2a+b=2g-i \\ b+2c=i}}p_1^*w^{[a]}\cdot p_2^*w^{[c]}\cdot \sum\limits_{d+e+f=b}(-1)^{d+f}p_1^*w^{[d]}\cdot \mu^*w^{[e]}\cdot p_2^*w^{[f]}
$$
where $p_1,p_2,\mu:J\times J\to J$ are the two projections and the multiplication map,
respectively.

We will use the following key property of the cycles $\pi_i$, established in \cite{Suh}:

\bpr
The endomorphism $\H^j(J_{k_\s},\Z/n)\to \H^j(J_{k_\s},\Z/n)$ of cohomology of $J_{k_\s}$ induced by the correspondence $\pi_i$ is zero for $j\neq i$, and is the identity map for $j=i$.
\epr
{\em Proof.}
This is proved in \cite[Theorem 4.2.3]{Suh} in the case $k_\s=\mathbb{C}$. The proof goes through over an arbitrary separably closed field: its only geometric input is Poincar\'e's formula $r!\cdot w^{[r]}=(w^{[1]})^{r}$ in cohomology of $J$, see \cite[Remark 2A13]{kleiman} for a proof of this over arbitrary fields. \hfill $\Box$

\medskip

Let us now upgrade the action of $\pi_i$ from individual cohomology groups to the complex $\RGamma_{\et}(J_{k_\s},\Z/n)$, as an object of $D(\Gamma,\Z/n)$.

\ble\label{correspondences-on-etalecoh}
For a smooth proper geometrically integral variety $Y$ of dimension $d$ over the field $k$, a cycle $c\in Z^i(Y\times_k Y)$ gives rise to a morphism in $D(\Gamma,\Z/n)$
$$
[c]:\RGamma_{\et}(Y_{k_\s},\Z/n)\to \RGamma_{\et}(Y_{k_\s},\Z/n)(i-d)[2i-2d]
$$
such that the induced maps on cohomology $\H^r_{\et}(Y_{k_\s},\Z/n)\to \H^{r+2i-2d}_{\et}(Y_{k_\s},\Z/n(i-d))$ coincide with the usual action of $c$ via correspondences on individual cohomology groups.
\ele 
{\em Proof.} The two projections $p_1,p_2:Y\times Y\to Y$ define maps in $D(\Gamma,\Z/n)$: $$p_1^*:\RGamma_{\et}(Y_{k_\s},\Z/n)\to \RGamma_{\et}((Y\times Y)_{k_\s},\Z/n)$$ and $$p_{2*}:\RGamma_{\et}((Y\times Y)_{k_\s},\Z/n)\to \RGamma_{\et}(Y_{k_{\s}},\Z/n)(-d)[-2d],$$ where $(-d)$ refers to the twist by the $(-d)$-th power of the cyclotomic character $\Gamma\to \Aut_{\Z/n}(\mu_n)$.

There are cycle classes $\cl(c)\in \H^{2i}_{\et}(Y\times Y,(\Z/n)(i))$ 
in absolute \'etale cohomology of $Y\times Y$, as defined in \cite[Cycle, 2.2.10]{sga45}.
The absolute \'etale cohomology complex $\RGamma_{\et}(Y\times Y,(\Z/n)(i))$ can be identified with the continuous group cohomology complex 
$\RGamma_{\cont}(\Gamma, \RGamma_\et((Y\times Y)_{k_\s},(\Z/n)(i)))$ with coefficients in geometric \'etale cohomology, hence $\cl(c)$ corresponds to a map $\cl(c):(\Z/n)[-2i]\to \RGamma_\et((Y\times Y)_{k_\s},(\Z/n)(i))$ in the derived category $D(\Gamma,\Z/n)$.

We can now define the endomorphism of the complex $\RGamma_{\et}(Y_{k_\s},\Z/n)$ induced by $c$ as the composition
\begin{multline*}
\RGamma_{\et}(Y_{k_\s},\Z/n)\xrightarrow{p_1^*} \RGamma_{\et}((Y\times Y)_{k_\s},\Z/n)\xrightarrow{\cup \cl(c)}\\ \RGamma_{\et}((Y\times Y)_{k_\s},\Z/n)(i)[2i]\xrightarrow{p_{2*}}\RGamma_{\et}(Y_{k_\s},\Z/n)(i-d)[2i-2d]
\end{multline*}
Here the middle map $\cup \cl(c)$ denotes the composition
\begin{multline*}
\RGamma_{\et}((Y\times Y)_{k_\s},\Z/n)\xrightarrow{\id\otimes \cl(c)}\RGamma_{\et}((Y\times Y)_{k_\s},\Z/n)^{\otimes 2}(i)[2i]\xrightarrow{}\\
\RGamma_{\et}((Y\times Y)_{k_\s},\Z/n)(i)[2i]
\end{multline*}
where the second map is the cup-product on the level of cohomology complexes. 
\hfill $\Box$

\medskip

We can now construct the desired quasi-isomorphism (\ref{jacobian-dec-formula}). It suffices to produce, for each $i\geq 0$, a map $\RGamma_{\et}(J_{k_\s},\Z/n)\to \H^i_{\et}(J_{k_\s},\Z/n)[-i]$ in $D(\Gamma,\Z/n)$ that induces the identity map on the $i$-th cohomology. 
By Lemma \ref{correspondences-on-etalecoh}, the cycle $\pi_i$ defines a map $[\pi_i]:\RGamma_{\et}(J_{k_\s},\Z/n)\to \RGamma_{\et}(J_{k_\s},\Z/n)$ in $D(\Gamma,\Z/n)$ inducing $0$ on $\H^j$ for $j\neq i$, and the identity map on $\H^i$. The endomorphism $[\pi_i]$ fits into the following map of distinguished triangles
$$
\begin{tikzcd}
\tau^{\leq i}\RGamma_{\et}(J_{k_\s},\Z/n)\arrow[r]\arrow[d,"{\tau^{\leq i}[\pi_i]}"] & \RGamma_{\et}(J_{k_\s},\Z/n) \arrow[d,"{[\pi_i]}"]\arrow[r] & \tau^{\geq i+1}\RGamma_{\et}(J_{k_\s},\Z/n)\arrow[d,"{\tau^{\geq i+1}[\pi_i]}"] \\
\tau^{\leq i}\RGamma_{\et}(J_{k_\s},\Z/n)\arrow[r] & \RGamma_{\et}(J_{k_\s},\Z/n) \arrow[r] & \tau^{\geq i+1}\RGamma_{\et}(J_{k_\s},\Z/n) 
\end{tikzcd}
$$
The map $\tau^{\geq i+1}[\pi_i]$ induces the zero map on all cohomology groups, 
but a priori it may be a non-zero map in $D(\Gamma,\Z/n)$. Since the complex $\tau^{\geq i+1}\RGamma_{\et}(J_{k_\s},\Z/n)$ has non-zero cohomology only in the range $[i+1,2g]$, the $(2g-i)$-th power of the endomorphism $\tau^{\geq i+1}[\pi_i]$ is nonetheless equal to zero in $D(\Gamma,\Z/n)$,
see Lemma \ref{nilpotent} below.

Therefore $[\pi_i]^{2g-i}$ gives rise to a map $\RGamma_{\et}(J_{k_\s},\Z/n)\to \tau^{\leq i}\RGamma_{\et}(J_{k_\s},\Z/n)$ inducing the identity map on $i$-th cohomology. Composing it with the natural map 
$$\tau^{\leq i}\RGamma_{\et}(J_{k_\s},\Z/n)\to \tau^{[i,i]}\RGamma_{\et}(J_{k_\s},\Z/n)\cong \H^i_{\et}(J_{k_\s},\Z/n)[-i]$$ 
we obtain the desired map $\RGamma_{\et}(J_{k_\s},\Z/n)\to \H^i_{\et}(J_{k_\s},\Z/n)[-i]$ in $D(\Gamma,\Z/n)$ that induces the identity map on $\H^i$. Summing up these maps over all $i\in [0,2g]$ we obtain the quasi-isomorphism (\ref{jacobian-dec-formula}).

Applying the derived functor of $\Gamma$-invariants to (\ref{jacobian-dec-formula}) gives the isomorphisms (\ref{jacobian-dec-oncoh-formula}).\hfill $\Box$

\medskip

\ble \label{nilpotent}
Let $A$ be an abelian category and let $D^b(A)$ be the bounded derived category of $A$.
Let $C$ be an object of $D^b(A)$
such that $\H^i(C)=0$ for $i\notin[a,a+n]$, for some integers $a$ and $n\geq 0$. 
Let $f\colon C\to C$ be an endomorphism in $D^b(A)$ such that $\H^i(f)=0$ for every $i$.
Then we have $f^{n+1}=0$ in $\Hom_{D^b(A)}(C,C)$. 
\ele
{\em Proof.} Without loss of generality we can assume $a=0$.
Proceed by induction on $n$. If $n=0$, then $C$  is isomorphic in $D^b(A)$
to an object represented by the one-term complex $\H^0(C)$ in degree $0$.
Under this isomorphism,  $f$ corresponds to an endomorphism 
induced by $\H^0(f)\in\Hom_A(\H^0(C),\H^0(C))$, but $\H^0(f)=0$.

Now suppose that $f^n$ is zero on $\tau^{<n}C$. 
The composition $C\xrightarrow{f} C \to \H^n(C)[-n]$
is zero, because it is equal to the composition $C\to \H^n(C)[-n]\to \H^n(C)[-n]$,
where the last map is $\H^n(f)=0$. Given the distinguished triangle 
$$\tau^{<n}C\to C\to\H^n(C)[-n],$$ this implies that $f\colon C\to C$
factors as $C\xrightarrow{g}\tau^{<n}C\to C$ for some map $g$. 
Then $f^{n+1}$ is the composition $C\xrightarrow{g}\tau^{<n}C\to C\xrightarrow{f^{n}}C$. 
The truncation $\tau^{<n}C$ is functorial in $C$, so the composition of the last two maps is equal to $\tau^{<n}C\xrightarrow{f^{n}}\tau^{<n}C\to C$. By the induction assumption 
$f^n$ is zero on $\tau^{<n}C$, so the iterate $f^{n+1}$ is zero on $C$. \hfill$\Box$

\section{Abelian surfaces over $\Q$ with $\delta_2^{0,2}\neq 0$}

In this section we give examples of abelian varieties $A$ over $\Q$ such that the differential
$\delta_2^{0,2}$ of the spectral sequence
\begin{equation}
E_{2}^{p,q}= \H^p(k,\H_\et^q(A_{k_\s},\Z_2(1)))\Rightarrow\H_\et^{p+q}(A,\Z_2(1)),
\label{ss twist}
\end{equation}
and the same differential of the spectral sequence with coefficients $\Z/2$, are both non-zero. 

Let $A$ be an abelian variety over a field $k$. We write $A^\vee$ for the dual abelian variety. We
denote by $\beta_2: \H^1(k, A)\to \H^2(k, A[2])$ the connecting homomorphism induced by the short exact sequence of $\Ga$-modules
$$0\to A[2](k_{\s})\to A(k_{\s})\xrightarrow{2}A(k_{\s})\to 0.$$

We use the following crucial proposition. 

\bpr \label{sunday}
Let $A$ be an abelian variety over a field $k$ of characteristic zero with a polarisation 
$$\lambda\in\NS(A_{k_\s})^\Ga=\Hom_k(A,A^\vee)^{\rm sym}.$$
Let $c'_\lambda\in\H^1(k,A^\vee)$ be the image of $\lambda$ under the connecting map
of the exact sequence of $\Ga$-modules
\begin{equation}
0\lra A^\vee(k_\s) \lra \Pic(A_{k_\s}) \lra \NS(A_{k_\s}) \lra 0.
\label{g1}
\end{equation}
Let ${\rm c}_1(\lambda)\in\H^2_\et(A_{k_\s},\Z_2(1))^\Ga$ be the first Chern class of
$\lambda$. 

Then the image of $\delta_2^{0,2}({\rm c}_1(\lambda))\in \H^2(k,\H^1_{\et}(A_{k_{\s}},\Z_2(1)))$ in $\H^2(k, \H^1_{\et}(A_{k_{\s}},\Z/2))$ is equal to $\beta_2(c'_{\lambda})\in \H^2(k,A^\vee[2])$. In particular, 
$\delta_2^{0,2}({\rm c}_1(\lambda))$ is divisible by $2$ in
$\H^2(k,\H^1_{\et}(A_{k_{\s}},\Z_2(1)))$ if and only if $c'_{\lambda}$  is
divisible by $2$ in $\H^1(k, A^{\vee})$.
\epr
{\em Proof.} The antipodal involution $[-1]\colon A\to A$  
induces an action of $\Z/2$ on $\Pic(A_{k_\s})$
which turns (\ref{g1}) into
an exact sequence of abelian groups with an action of $\Z/2$. The induced action on 
$\NS(A_{k_\s})$ is trivial, cf.~\cite[Remark 12.13]{MilAV}. The involution $[-1]_A$
induces the involution $[-1]_{A^\vee}$ on $A^\vee$.
Recall that for an abelian group $M$ with an action of $\Z/2=\langle\sigma\rangle$ 
we have $\H^1(\Z/2,M)\cong M^{\si=-1}/(1-\si)M$. 
Since $A^\vee(k_\s)$ is 2-divisible, the endomorphism $1-[-1]_{A^{\vee}}=[2]_{A^{\vee}}$ is surjective, and we obtain $\H^1(\Z/2,A^\vee(k_\s))=0$. 
Thus the long exact sequence of cohomology for the group $\Z/2$ gives an exact sequence of $\Gamma$-modules
\begin{equation}
0 \lra A^\vee[2] \lra \Pic(A_{k_\s})^{[-1]^*} \lra \NS(A_{k_\s}) \lra 0.
\label{g2}
\end{equation}
Let $c_\lambda\in \H^1(k,A^\vee[2])$ be the image of 
$\lambda$ under the connecting map 
of $(\ref{g2})$. It is clear that $c'_\lambda$ is the image of $c_\lambda$ in $\H^1(k,A^\vee)$.

Let $T_2(A^\vee)$ be the 2-adic Tate module, and let $e_2\colon A[2]\times A^\vee[2]\to\Z/2$ be the Weil pairing.
The short exact sequences (\ref{g2}) and (\ref{alpha}) are compatible,
so that there is a commutative diagram (see diagram (16) in \cite{PR}):
$$\xymatrix{
0 \ar[r]& A^\vee[2] \ar[r]\ar[d]^\id& \Pic(A_{k_\s})^{[-1]^*} \ar[r]\ar[d]& \NS(A_{k_\s}) \ar[r]\ar[d]& 0\\
0 \ar[r]& A^\vee[2] \ar[r]&S^2(A^\vee[2])\ar[r]&\wedge^2(A^\vee[2])\ar[r]&0}$$
Furthermore, the right-hand vertical map sends $\lambda$ to the element of $\wedge^2(A^\vee[2])$
corresponding to 
$$e_2(x,\lambda(y))\in \Hom(\wedge^2A[2],\Z/2).$$
This map factors as the first Chern class map ${\rm c}_1$ 
$$\NS(A_{k_\s})\cong\Hom(A_{k_\s},A^\vee_{k_\s})^{\rm sym}\stackrel{{\rm c}_1}\lra\Hom(T_2(A),T_2(A^\vee))^{\rm sym}\cong\wedge^2 T_2(A^\vee)(-1)
\cong\H^2_\et(A_{k_\s},\Z_2(1))$$
followed by the reduction modulo 2 map
$\H^2_\et(A_{k_\s},\Z_2(1))\to\H^2_\et(A_{k_\s},\Z/2)$, see \cite[Lemma 2.6]{OSZ}.
Theorem \ref{p1} gives 
\begin{equation}\label{main formula on clambda}
\delta_2^{0,2}({\rm c}_1(\lambda))={\rm Bock}_{T_2(A^\vee)}(c_\lambda),
\end{equation}
where ${\rm Bock}_{T_2(A^\vee)}$ is the connecting map
attached to the exact sequence
$$0\to T_2(A^\vee)\xrightarrow{2} T_2(A^\vee)\to A^\vee[2]\to 0.$$
This sequence fits into the following commutative diagram with exact rows:
$$\xymatrix{0\ar[r]&A^\vee[2]\ar[r]& A^\vee\ar[r]^{2} &A^\vee\ar[r]& 0\\
0\ar[r]& A^\vee[2]\ar[r]\ar[u]^\cong& A^\vee[4]\ar[r]\ar[u] &A^\vee[2]\ar[r]\ar[u]& 0\\
0\ar[r]& T_2(A^\vee)\ar[r]^{2}\ar[u]& T_2(A^\vee)\ar[r]\ar[u]& A^\vee[2]\ar[r]\ar[u]_\cong& 0}$$
This implies a relation between connecting homomorphisms induced by the first and last rows: the image of $\Bock_{T_2(A^{\vee})}(c_{\lambda})\in \H^2(k, T_2A^{\vee})$ in $\H^2(k, A^{\vee}[2])$ equals $\beta_2(c'_{\lambda})$, where $\beta_2:\H^1(k, A^{\vee})\to \H^2(k, A^{\vee}[2])$ is the connecting homomorphism induced by the top exact row, as defined above. Combining this with (\ref{main formula on clambda}) gives the desired equality. 

To prove the final assertion of the proposition, note that the short exact sequence of Galois modules $$0\to \H^1_{\et}(A_{k_{\s}},\Z_2(1))\xrightarrow{2}\H^1_{\et}(A_{k_{\s}},\Z_2(1))\to \H^1_{\et}(A_{k_{\s}},\Z/2)\to 0$$
shows that $\delta^{0,2}_2(c_1(\lambda))\in \H^2(k,\H^1_{\et}(A_{k_{\s}},\Z_2(1)))$ is divisible by $2$ if and only if its image in $\H^2(k, \H^1_{\et}(A_{k_{\s}},\Z/2))\cong \H^2(k, A^{\vee}[2])$ is zero. Similarly, the short exact sequence of Galois modules
$$0\to A^{\vee}[2]\to A^{\vee}(k_{\s})\xrightarrow{2}A^{\vee}(k_{\s})\to 0$$
gives that $c_{\lambda}'\in \H^1(k,A^{\vee})$ is divisible by $2$ if and only if its image $\beta_2(c'_{\lambda})\in \H^2(k,A^{\vee}[2])$ under the connecting homomorphism $\beta_2$ is zero. This implies the equivalence of $2$-divisibility of $\delta_{2}^{0,2}(c_1(\lambda))$ and $c'_{\lambda}$.
\hfill $\Box$

\medskip

Proposition \ref{sunday} demonstrates that the Hochschild--Serre spectral sequence for $A$ (both with $\Z_2$ and $\Z/2$ coefficients) does not degenerate if the class $c'_\lambda\in \H^1(k,A^{\vee})$ is not divisible by $2$. In the course of the above proof we also introduced the class $c_{\lambda}\in \H^1(k,A^{\vee}[2])$ lifting $c_{\lambda}'$.
Classes $c_{\lambda}$ and $c_{\lambda}'$ have an explicit interpretation when $A$ is the Jacobian of a curve, which
is particularly convenient when $C$ is hyperelliptic.

Let $C$ be a smooth, projective, geometrically integral curve of genus $g$ over a field $k$
of characteristic not equal to 2.
Let $J:={\bf Pic}^0_{C/k}$ be the Jacobian of $C$ with its canonical principal
polarisation $\lambda$.
By \cite[Theorem~3.9]{PR}, $c_\lambda$ is the class of the torsor of theta characteristics,
which is the closed subvariety of 
${\bf Pic}_{C/k}^{g-1}$ given by $2x=K_C$, where $K_C$ is the canonical class.
Thus $c'_\lambda$ is the class of ${\bf Pic}_{C/k}^{g-1}$ in $\H^1(k,J)$.

If $C$ is a hyperelliptic curve of odd genus or with
a rational Weierstrass point, then $c_\lambda=0$, see \cite[Proposition~3.11]{PR}.
If $C$ is a hyperelliptic curve of even genus,
then we have an isomorphism ${\bf Pic}_{C/k}^{g-1}\cong{\bf Pic}_{C/k}^{1}$
given by subtracting $\frac{g-2}{2} H_C$, 
where $H_C$ is the hyperelliptic divisor class. Then $c_\lambda$ is the class of the torsor
for $J[2]$ in ${\bf Pic}_{C/k}^{1}$ given by $2x=H_C$, 
and $c'_\lambda=[{\bf Pic}_{C/k}^{1}]$. Explicitly, the $k_\s$-points of the torsor
of theta characteristics are Galois-equivariantly identified with the subsets of 
the set of Weierstrass points of $C$ of odd cardinality, after identifying each subset with its complement.

\medskip

Now let $k$ be a number field. In this case $c_\lambda$ lies in the 
2-Selmer group $\Sel_2(J)$,
see \cite[Corollary~2]{PS}. Thus $c'_\lambda$ lies in the Tate--Shafarevich group $\Sha(J)$.

We have the Cassels--Tate pairing 
$\Sha(J)\times \Sha(J^\vee)\to\Q/\Z$, see \cite[Ch.~I, \S 6]{ADT}.
Identifying $J$ with its dual abelian variety $J^\vee$ via $\lambda$, we obtain a pairing
$$\langle x,y\rangle \colon \Sha(J)\times \Sha(J)\to\Q/\Z,$$
which we shall also call the Cassels--Tate pairing.
By \cite[Theorem~5]{PS} we have $\langle x,x+c'_\lambda\rangle=0$ for any $x\in \Sha(J)$.

For a finite $\Ga$-module $M$ one defines
$\Sha^1(k,M)$ as the subgroup of $\H^1(k,M)$
consisting of the classes that restrict to zero in $\H^1(k_v,M)$,
for all places $v$ of $k$. It is clear that $\Sha^1(k,J[2])\subset\Sel_2(J)$.

Over a number field $k$, we can test the divisibility of $c'_\lambda$ in $\H^1(k,J)$
in terms of the Cassels--Tate pairing.
Let us denote by $i$ the natural map $\H^1(k,J[2])\to\H^1(k,J)$. 

Recall that a curve $C$ over a number field $k$ is called {\em odd} if
$\langle c'_\lambda,c'_\lambda\rangle\neq 0$. 

\ble \label{monday}
The class $c'_\lambda$ is not divisible by $2$ in $\H^1(k,J)$ if and only if 
$i(\Sha^1(k,J[2]))$ pairs non-trivially with $c'_\lambda$ under the Cassels--Tate pairing.
For example, this holds when $C$ is odd and $c_\lambda\in\Sha^1(k,J[2])$.
\ele
{\em Proof.} This is a particular case of \cite[Theorem~4]{Cre13}, which is
based on the non-degeneracy of the Poitou--Tate pairing \cite[Ch.~I, \S 4]{ADT}
$$\Sha^1(k,J[2])\times \Sha^2(k,J^\vee[2])\to\Z/2$$
and its compatibility with the Cassels--Tate pairing
$ \Sha(J)\times \Sha(J^\vee)\to\Q/\Z$. For the second claim we consider the
pairing of $c'_\lambda=i(c_\lambda)$ with itself. \hfill $\Box$

\brem{\rm \label{poonenrains}
Suppose that $C$ is a hyperelliptic curve with equation $y^2=f(x)$, 
where $f(x)\in k[x]$ is a separable polynomial. If $f(x)$ has a root in $k_v$
for every place $v$ where $C$ has bad reduction
or the residue characteristic of $k_v$ is $2$, then $c_\lambda\in\Sha^1(k,J[2])$, see 
\cite[Proposition 3.6]{PR}.
}\erem

Manin pointed out that the Cassels--Tate pairing can be interpreted as a particular case of 
what is now called the Brauer--Manin pairing. Let us recall this interpretation.
Let $\Br(C)$ be the Brauer group of $C$ and let $\Br_0(C)=\Im[\Br(k)\to\Br(C)]$.
Since $\Br(C_{k_\s})=0$, the spectral sequence 
$\H^p(k,\H^q_\et(C_{k_\s},\G_m))\Rightarrow\H^{p+q}_\et(C,\G_m)$
gives an exact sequence
\begin{equation}
0\to \Br_0(C)\to\Br(C)\to\H^1(k,\Pic(C_{k_\s}))\to 0, \label{night}
\end{equation}
where we used that $\H^3(k,k_\s^\times)=0$ which holds since $k$ is a number field.
We have an exact sequence of $\Gamma$-modules
$$0\to J(k_\s)\to \Pic(C_{k_\s})\xrightarrow{\deg}\Z\to 0.$$
It induces an isomorphism of
$\H^1(k,\Pic(C_{k_\s}))$ with the quotient of  $\H^1(k,J)$ by the cyclic subgroup 
generated by $[{\bf Pic}_{C/k}^{1}]$. 
For an element $x\in\H^1(k,J)$
we say that a Brauer class ${\mathcal A}\in \Br(C)$ is {\it associated} to $x$ if the images
of $x$ and ${\mathcal A}$ in $\H^1(k,\Pic(C_{k_\s}))$ are equal.
If $\mathcal{A}$ is associated to $x$, then there is a natural choice of an extension of $\mathcal{A}$ to a Brauer class on ${\bf Pic}_{C/k}^{1}$.

\ble\label{Cr}
Let $k$ be a number field and let $C$ be an everywhere locally soluble 
hyperelliptic curve of even genus over $k$. If the class in $\Br(C)$ 
associated to some element of $i(\Sha^1(k,J[2]))\subset\H^1(k,J)$
obstructs the Hasse principle on $C$, then $c'_\lambda$ is not divisible by $2$ in $\H^1(k,J)$.
\ele
{\em Proof.} Since $C$ is hyperelliptic of even genus, 
we have ${\bf Pic}_{C/k}^{g-1}\cong {\bf Pic}_{C/k}^{1}$.
For each place of $k$ choose a local point $P_v\in C(k_v)\subset {\bf Pic}_{C/k}^{1}(k_v)$.
By a theorem of Manin \cite[Theorem~6.2.3]{S01}, we have
$$\sum_v\inv_v({\mathcal A}(P_v))=\langle x,c'_\lambda\rangle\in\frac{1}{2}\Z/\Z,$$
for any $x\in\Sha(J)$ and  ${\mathcal A}\in \Br(C)$ associated to $x$, 
using that $c'_\lambda=[{\bf Pic}_{C/k}^{1}]\in\Sha(J)$. Taking $x\in i(\Sha^1(k,J[2]))$
we conclude by Lemma \ref{monday}. \hfill $\Box$

\bthe\label{nonsplit abelian example}
Let $J$ be the Jacobian of the hyperelliptic curve of genus $2$ over $\Q$ given by

{\rm (a)} $y^2=3(x^2+1)(x^2+17)(x^2-17)$, or

{\rm (b)} $y^2 = -(x^6+x+6).$

\noindent Then the differential $\delta_2^{0,2}$ of the spectral sequence $(\ref{Leray2})$
and the same differential of a similar spectral sequence with coefficients $\Z/2$ are both non-zero.
\ethe
{\em Proof.} (a) By Creutz \cite[p.~941]{Cre13} and Creutz and Viray \cite[Theorem~6.7]{CV15},
this curve satisfies the assumption of Lemma \ref{Cr}. Here is a sketch of this calculation.

Let $C$ be the hyperelliptic curve over a field $k$ of characteristic different from 2
given by $y^2=cf(x)$ where $c\in k^\times$ and $f(x)\in k[x]$ is a separable monic
polynomial of even degree. Let $J$ be the Jacobian of $C$. Let $L=k[x]/(f(x))$ and let
$\theta\in L$ be the image of $x$. The well-known identification of the group $k$-scheme $J[2]$
with the quotient of the kernel of the norm map $R_{L/k}(\mu_2)\to \mu_2$ by the
image of $\mu_2$ gives a natural inclusion
$$(L^\times/k^\times L^{\times 2})_1\subset\H^1(k,J[2])/\langle c_\lambda\rangle,$$
where the subscript $1$ denotes the subgroup of elements with norm 
$1\in k^\times/k^{\times 2}$. To $l\in L^\times$ we associate the class of the quaternion algebra
$(l,x-\theta)\in\Br(L(x))$, which only depends on the image of $l$ in $L^\times/L^{\times 2}$.
The corestriction
$${\rm cores}_{L(x)/k(x)}\big((l,x-\theta)\big)\in \Br(k(\P^1_k))$$
is unramified away from the ramification locus of $C\to \P^1_k$ and the point at infinity.
If the norm of $l$ is a square in $k^\times$, then it is unramified at infinity.
It follows that the pullback of this element to
$\Br(k(C))$ is unramified, and so belongs to the subgroup $\Br(C)\subset \Br(k(C))$.
We denote the resulting element by ${\mathcal A}_l\in\Br(C)$.
Finally, multiplying $l$ by $s\in k^\times$ gives ${\mathcal A}_{ls}={\mathcal A}_l+(s,c)$,
so this does not change the image of ${\mathcal A}_l$ in $\Br(C)/\Br_0(C)$. 
An important property
of this construction is that the map sending $l$ to ${\mathcal A}_l$ coincides with the composition 
$$\H^1(k,J[2])/\langle c_\lambda\rangle\stackrel{i}\lra \H^1(k,J)/\langle c'_\lambda\rangle\lra \Br(C)/\Br_0(C),$$
where the second arrow comes from (\ref{night}).

Now let $k=\Q$ and $L=\Q(\sqrt{-1})\oplus\Q(\sqrt{-17})\oplus\Q(\sqrt{17})$.
The 0-dimensional scheme $\Spec(L)$ over $\Q$
is everywhere locally soluble, hence the curve $C$ given by
$$y^2=c(x^2+1)(x^2+17)(x^2-17)$$ 
is everywhere locally soluble for any $c\in \Q^\times$ and we have $c_\lambda\in\Sha^1(\Q,J[2])$.
Take $l=(1,1,-1)\in L$.
One checks that $l$ gives rise to an everywhere locally trivial element of 
$(L^\times/\Q^\times L^{\times 2})_1$, and thus to an element of 
$\Sha^1(\Q,J[2])/\langle c_\lambda\rangle$.
We have ${\mathcal A}_l=(-1,x^2-17)$. A local computation \cite[Lemma~6.8]{CV15} shows that if
the number of odd prime factors $p$ of $c>0$ such that neither $17$ nor $-1$
is a square modulo $p$, is odd, then 
$\sum_v\inv_v({\mathcal A}_l(P_v))=1/2$. For example, one can take $c=3$.
Thus Lemma \ref{Cr} can be applied, so that $c'_\lambda$ is not divisible by $2$ in 
$\H^1(\Q,J)$, hence $\delta^{0,2}_2$ is non-zero by Proposition \ref{sunday}. 

\smallskip

(b) This curve was also communicated to us by Creutz (personal communication).
It is the quadratic twist by $-1$ of the curve in \cite[Example 3.12]{PR}. 
The arguments in {\em loc.~cit.}~show that for the primes $p$ where $C$ has bad reduction, 
and also for $p=2$, the polynomial $f(x)$ has a root in~$\Q_p$.
In view of Lemma \ref{monday} and Remark \ref{poonenrains}
it remains to show that $C$ is odd. By \cite[Corollary 12]{PS} this happens if and only if 
the number of places $v$ such that $C$ does not have a $\Q_v$-rational divisor of degree 1, is odd.
Such places are necessarily real or places of bad reduction, so we only need to consider the real place.
We check that $C(\R)=\emptyset$, so $C$ has no real divisor of degree 1, and hence is odd.
\hfill $\Box$

\section{The Brauer group of a torus} 

Let $T$ be a torus over a field $k$ of characteristic exponent $p$.
Let $\widehat{T}=\Hom_{k_\s-{\rm gps}}(T_{k_\s},\G_m)$ be the $\Ga$-module of characters of $T$.
The Brauer group $\Br(T)=\H^2_\et(T,\G_m)$ is computed by the spectral sequence
\begin{equation}
E_{2}^{p,q} = \H^p(k, \H^q_{\et}(T_{k_\s}, \G_{m})) \Rightarrow \H^{p+q}_{\et}(T, \G_{m}).
\label{Leray-torus}
\end{equation}
We denote by $d^{p,q}_r$ the differential on the $r$th page emanating from the $(p,q)$-entry.

Since $T_{k_\s}$ is a dense open subscheme of
$\A^n_{k_\s}$, we have $\H^1(T_{k_\s},\G_{m})\cong\Pic(T_{k_\s})=0$.
Thus $d^{0,2}_2=0$, and the first interesting differential
is $d^{\,0,2}_3\colon\Br(T_{k_\s})^\Ga\to \H^3(k,k_\s[T]^\times)$.

The origin $e\in T(k)$
of the group law on $T$ gives a section of the structure morphism $T\to\Spec(k)$,
hence $\H^r(k,\G_m)\to\H^r_\et(T,\G_m)$ is injective for all $r\geq 0$.
Likewise, the natural map $\Br(k)\to\Br(T)$ is injective.
Let $\Br_e(T) =\Ker [\Br(T) \to \Br(k)]$ be the kernel of specialisation at $e$.

We have an isomorphism of $\Ga$-modules
$\H^0_\et(T_{k_\s},\G_{m})=k_\s[T]^\times\cong k_\s^\times \oplus \widehat{T}$ (cf. \cite[Proposition 3]{Ros}).
The spectral sequence (\ref{Leray-torus}) thus gives rise to an exact sequence
\begin{equation} 
0 \to \H^2(k,\widehat{T}) \to \Br_e(T) \to \Br(T_{k_\s})^\Ga \xrightarrow{\bar d_3^{\,0,2}} \H^3(k, \widehat{T}). \label{ias1}
\end{equation}
To compute $\Br(T)$ one needs to describe the map 
$\bar d_3^{\, 0,2}\colon \Br(T_{k_\s})^\Ga \to \H^3(k, \widehat{T})$, which is the composition of the differential $d_3^{0,2}$ with the map induced by the projection $k_{\s}[T]^{\times}\to \widehat{T}$.

Let $\ell$ be a prime number not equal to $p$. For an abelian group $A$ we write $A(p')$
for the subgroup of $A$ consisting of the elements of finite order not divisible by $p$.

Since $\Pic(T_{k_\s})=0$, the Kummer sequence gives rise to an isomorphism
$$\kappa\colon\H^2_\et(T_{k_\s},\mu_{\ell^n})\tilde\lra \Br(T_{k_\s})[\ell^n].$$ 
Using the isomorphism of $\Ga$-modules 
$\H^0_\et(T_{k_\s},\G_{m})\cong k_\s^\times \oplus \widehat{T}$,
we also deduce from the Kummer sequence a natural isomorphism 
$$\widehat T/\ell^n\tilde\lra \H^1_\et(T_{k_\s},\mu_{\ell^n}).$$
We note that $\H^2_\et(T_{k_\s},\Z/\ell^n)\cong\wedge^2\H^1_\et(T_{k_\s},\Z/\ell^n)$
and thus the multiplication by $m$ map $[m]\colon T\to T$ acts on $\H^2_\et(T_{k_\s},\mu_{\ell^n})$,
and hence on $\Br(T_{k_\s})(p')$, as $m^2$. On the other hand, $[m]$ acts on $\widehat{T}$
as $m$. Taking $m=-1$ we see that $2\bar d_3^{\, 0,2}=0$, so that $\bar d_3^{\, 0,2}$
sends the elements of $\Br(T_{k_\s})(p')$ of odd order to zero.
Then it follows from (\ref{ias1}) that every element of odd order in
$\Br(T_{k_\s})(p')^\Ga$ lifts to $\Br(T)$.

The question of an explicit description of the map $\bar d_3^{\, 0,2}$ 
on the 2-primary torsion subgroup of $\Br(T_{k_\s})^\Ga$
was asked at the top of \cite[p.~220]{CTS21}. We now give such a description.

\bthe \label{tori: main gm coeffs}
Let $T$ be a torus over a field $k$ of characteristic different from $2$. Let $\bar d^{\,0,2}_3$ be the composition
$$\Br(T_{k_\s})^\Ga\stackrel{d^{\,0,2}_3}\lra \H^3(k,k_\s[T]^\times)\to
\H^3(k,\widehat T),$$
where the last map is induced by the projection  
$k_\s[T]^\times\cong k_\s^\times\times\widehat T\to \widehat T$. 
The restriction of $\bar d^{\,0,2}_3$ to 
the $2^n$-torsion subgroup is the composition
$$\Br(T_{k_\s})[2^n]^\Ga\xrightarrow[\simeq]{\kappa^{-1}}\wedge^2(\widehat T/2^n)(-1)^\Ga 
\stackrel{\rm Bock}\lra \H^1(k,\wedge^2(\widehat T/2))
\stackrel{\alpha[1]}\lra\H^2(k,\widehat T/2)\stackrel{\rm Bock}\lra \H^3(k,\widehat T).$$
\ethe
{\em Proof.} We will show that the restriction of $d_3^{0,2}$ to the $2^n$-torsion subgroup can be read off from the Hochschild--Serre spectral sequence for the cohomology of $T$ with coefficients in $\mu_{2^n}$. The natural map $\tau^{\leq 2}\RGamma_{\et}(T_{k_{\s}},\mu_{2^n})\to \tau^{\leq 2}\RGamma_{\et}(T_{k_{\s}},\G_m)$ of truncations in degrees $\leq 2$ induces a commutative diagram relating the extensions between $\H^2$ and $\tau^{\leq 1}$ in these complexes:
\begin{equation}\label{torus: kummer diagram}
\begin{tikzcd}
\H^2_{\et}(T_{k_{\s}},\mu_{2^n})\arrow[r]\arrow[d] & \tau^{\leq 1}\RGamma(T_{k_{\s}},\mu_{2^n})[3]\arrow[d] \\
\H^2_{\et}(T_{k_{\s}},\G_m)\arrow[r] & \H^0_\et(T_{k_{\s}},\G_m)[3]
\end{tikzcd}
\end{equation}
where we used that $\H^1_\et(T_{k_{\s}},\G_m)=\Pic(T_{k_{\s}})$ vanishes. The differential $d_{3}^{0,2}$ is obtained by applying the functor $\H^0(k,-)$ to the bottom horizontal map in the diagram.

Applying the functor of Galois invariants to the left vertical map we get the inclusion of the $2^n$-torsion 
$\wedge^2(\hT/2^n)(-1)^{\Gamma}\cong \H^2_{\et}(T_{k_{\s}},\mu_{2^n})$ into the Galois invariants in the Brauer group of $T_{k_{\s}}$, so we are looking to compute the result of applying 
$\H^0(k,-)$ to the counter-clockwise composition in (\ref{torus: kummer diagram}), composed with the projection $\H^3(k,\H^0_\et(T_{k_{\s}},\G_m))\to \H^3(k,\hT)$.

From the Kummer sequence, we see that $\tau^{\leq 1}\RGamma_{\et}(T_{k_{\s}},\mu_{2^n})$ can be represented by the two-term complex $\H^0_\et(T_{k_{\s}},\G_m)\xrightarrow{2^n}\H^0_\et(T_{k_{\s}},\G_m)$ with the right vertical arrow in (\ref{torus: kummer diagram}) given by shift by $[3]$ of the projection onto the $0$th term of this complex. 
This map followed by the projection to $\hT$ sends $\H^0_\et(T_{k_\s},\mu_{2^n})$ to zero, so the composition factors through $\tau^{[1,1]}\RGamma_{\et}(T_{k_{\s}},\mu_{2^n})$.
Applying $\H^0(k,-)$ to this composition, we obtain a map that factors as follows:
$$\H^0(k, \tau^{\leq 1}\RGamma_\et(T_{k_{\s}},\mu_{2^n})[3])\to \H^2(k, \H^1_\et(T_{k_{\s}},\mu_{2^n}))\cong \H^2(k,\hT/2^n)\xrightarrow{-\Bock_{\hT}}\H^3(k, \hT).$$ 
Therefore the clockwise composition in (\ref{torus: kummer diagram}) evaluated on $\H^0(k,-)$, composed with the projection onto $\H^3(k,\hT)$, is given by
$$\H^2_\et(T_{k_{\s}},\mu_{2^n})^{\Gamma}\xrightarrow{\delta_{2}^{0,2}} \H^2(k, \H^1_\et(T_{k_{\s}},\mu_{2^n}))\xrightarrow{-\Bock_{\hT}}\H^3(k,\hT),$$
where $\delta^{0,2}_2$ is a differential of the spectral sequence
$$E^{p,q}_2=\H^p(k,\H^q_\et(T_{k_\s},\mu_{2^n}))\Rightarrow\H^{p+q}_\et(T,\mu_{2^n}).$$
By Theorem \ref{p2}, $\delta^{0,2}_2$ is obtained by applying $-\otimes_{\Z_2}^L\mu_{2^n}$ to
the terms of the 2-extension given by
the difference of the maps (\ref{genius2}) and (\ref{genius3}) with 
$M=\H^1_\et(T_{k_\s},\Z_2)$. The composition of two Bockstein maps is zero, 
so only (\ref{genius3}) contributes to $\bar d_3^{\,0,2}$, thus proving that $\bar d_3^{\,0,2}$
is the composition of the four maps in the theorem. \hfill$\Box$

\brem{\rm\label{tori: nonvanishing remark}
It would be interesting to construct a torus with a non-zero map $\bar d_3^{\, 0,2}$, or
prove that none exist. 
}
\erem

When the first version of this paper was completed, the authors became aware of the
following result of Julian Demeio \cite[Theorem 1.1]{Dem}.
Recall that a torus is called {\em quasi-trivial} if $\widehat T$ is a permutational $\Ga$-module,
that is, $\widehat T$ has a $\Ga$-stable $\Z$-basis.

\bco \label{demeio}
Let $k$ be a field of characteristic zero.
If $T$ is a quasi-trivial torus, or if $k$ is a local or global field,
then the natural map $\Br(T)\to\Br(T_{k_\s})^\Ga$ is surjective.
\eco
{\em Proof.} Let us show that $\bar d^{\,0,2}_3=0$.
The case of quasi-trivial torus is immediate
from Theorem \ref{tori: main gm coeffs} and Lemma \ref{permutational} (i). 
If $k$ is a $p$-adic field, then $\H^3(k,\widehat T)=0$ since $k$ has strict cohomological dimension 
2 \cite[Theorem 10.6]{Har20}, so there is nothing to prove.
The case $k_v=\R$ follows from 
Theorem \ref{tori: main gm coeffs} and Lemma \ref{permutational} (ii).
If $k$ is a number field, we have an isomorphism $\H^3(k,\widehat T)\tilde\lra\prod
\H^3(k_v,\widehat T)$, where the product is over the real places of $k$,
see \cite[Exercise 18.1]{Har20}, so this case follows from the case of local fields.
\hfill $\Box$

\section{Torsors} \label{Torsors}

Let $X$ be a $k$-torsor for a semiabelian variety $A$.
In this section we address the problem of computing the \'etale cohomology
groups $\H^i_\et(X,\Z/n)$, where $n$ is not divisible by ${\rm char}(k)$.

Translations by elements of $A(k_\s)$ act trivially on the groups
$\H^i_\et(X_{k_\s},\Z/n)$, for $i\geq 0$. Let us explain this.
The choice of a $k_\s$-point in $X$
gives an isomorphism $X_{k_\s}\simeq A_{k_\s}$.
By the canonical isomorphism $\H^i_\et(A_{k_\s},\Z/n)\cong\wedge^i \H^1_\et(A_{k_\s},\Z/n)$
it suffices to prove the triviality of action on $\H^1_\et(A_{k_\s},\Z/n)$.
Let $p_1,p_2,\mu:A\times A\to A$ be the two projections and the multiplication map,
respectively. We need to show that for any $a\in A(k_\s)$ the composition
$$\H^1_{\et}(A_{k_\s},\Z/n)\xrightarrow{\mu^*} \H^1_{\et}(A_{k_\s}\times A_{k_\s},\Z/n)\xrightarrow{(a\times\Id_{A})^*}\H^1_{\et}(A_{k_\s},\Z/n)$$
is the identity map. Using that $\H^1(k_\s,\Z/n)=0$, we see that it suffices to prove that
$\mu^*=p_1^*+p_2^*$, but this follows from the canonical isomorphism 
\begin{equation}
\H^1_\et(A_{k_\s},\Z/n)\cong\Hom(\pi_1^\et(A_{k_\s}),\Z/n)\cong\Hom(A[n],\Z/n), \label{isomo}
\end{equation}
see (\ref{25a}).

\medskip

Therefore, any choice of a trivialisation $X_{k_{\s}}\simeq A_{k_{\s}}$
of the torsor $X$ over $k_{\s}$ gives rise to the same $\Gamma$-equivariant isomorphism
$$\H^i_\et(X_{k_\s},\Z/n)\cong\H^i_\et(A_{k_\s},\Z/n), \quad i\geq 0.$$
But the complexes $\RGamma_{\et}(X_{k_{\s}},\Z/n)$ and $\RGamma_{\et}(A_{k_{\s}},\Z/n)$ need not be isomorphic as objects of $D(\Gamma,\Z/n)$: the truncation of the latter in degrees $[0,1]$ is the direct sum of its cohomology groups, but the truncation $\tau^{[0,1]}\RGamma_{\et}(X_{k_{\s}},\Z/n)$ need not be split, as we demonstrate in this section.

The Hochschild--Serre spectral sequence for $X$ has the form
\begin{equation}
E^{p,q}_2=\H^p(k,\H^q_\et(A_{k_\s},\Z/n))\cong\H^p(k,\H^q_\et(X_{k_\s},\Z/n))\Rightarrow\H^{p+q}_\et(X,\Z/n). \label{ssX}
\end{equation}
The question we address is the explicit form of the differentials
$$\delta^{p,1}_2\colon\H^p(k,\H^1_\et(A_{k_\s},\Z/n)) \to \H^{p+2}(k,\Z/n),$$
where $p\geq 0$. 
Each of these differentials is induced by the map 
in the derived category of discrete $\Z/n$-modules with a continuous action of $\Ga$
$$
    \delta_X:\H^1_{\et}(A_{k_{\s}},\Z/n)\cong \H^1_\et(X_{k_{\s}},\Z/n)\to \Z/n[2]
$$
arising from the complex $\tau^{[0,1]}\RGamma_{\et}(X_{k_{\s}},\Z/n)\in D(\Gamma,\Z/n)$. 

Recall that $\mathcal{A}(\Ga,\Z/n)$ is
the abelian category of discrete $\Z/n$-modules equipped with a continuous action of $\Ga$.
Writing $M=\H^1_{\et}(A_{k_{\s}},\Z/n)$, we can think of
$\delta_X$ as an element of
$\Ext^2_{\mathcal{A}(\Ga,\Z/n)}(M,\Z/n)$. Since
$\Hom_{\mathcal{A}(\Ga,\Z/n)}(M,\Z/n)=\Hom_{\Z/n}(M,\Z/n)^\Ga$, we have a spectral sequence
\begin{equation}
E^{p,q}_2=\H^p(k,\Ext^q_{\Z/n}(M,\Z/n))\Rightarrow\Ext^{p+q}_{\mathcal{A}(\Ga,\Z/n)}(M,\Z/n).
\label{ss_ext}
\end{equation}
In our case, $M$ is a free, hence projective $\Z/n$-module, thus $\Ext^q_{\Z/n}(M,\Z/n)=0$
for $q>0$. Now (\ref{ss_ext}) gives a natural isomorphism
\begin{equation}
\Ext^2_{\mathcal{A}(\Ga,\Z/n)}(M,\Z/n)\cong \H^2(k,\Hom_{\Z/n}(M,\Z/n))\cong
\H^2(k,A[n]), \label{iso_delta}
\end{equation}
where the second isomorphism is the $\Z/n$-linear dual of a natural isomorphism 
(\ref{isomo}) of $\Ga$-modules
$\Hom_{\Z/n}(A[n],\Z/n)\cong\H^1_\et(A_{k_\s},\Z/n)$ induced by
the map (\ref{pi1coh to etalecoh map}). Using isomorphisms (\ref{iso_delta}), we identify
$\delta_X$ with an element of $\H^2(k,A[n])$.

We will express the class $\delta_X$ in terms of the class of the torsor $X$. 
The exact sequence of $\Ga$-modules
\begin{equation}
0\to A[n]\to A\to A\to 0 \label{torsor}
\end{equation}
gives rise to the homomorphism of Galois cohomology groups
\begin{equation}\label{torsor bockstein}
\beta_n\colon \H^1(k,A)\to\H^2(k,A[n]).
\end{equation}

\bthe \label{1}
Let $k$ be a field and let $n$ be a positive integer not divisible by ${\rm char}(k)$.
Let $X$ be a $k$-torsor for a semiabelian variety $A$ over $k$ with class $[X]\in\H^1(k,A)$. The class $\delta_X\in \H^2(k,A[n])$ is equal to the image of $[X]$ under the map $\beta_n$. 
In particular, the differentials $\delta^{p,1}_2$ are given by cupping with the class 
$$\beta_n([X])\in \H^2(k, A[n])\cong \Ext^2_{\mathcal{A}(\Ga,\Z/n)}(\H^1(A_{k_{\s}},\Z/n),\Z/n).$$
\ethe

{\em Proof.}
Our goal is to compute the morphism 
$$\delta_X:\Hom_{\Z/n}(A[n],\Z/n)\cong 
\H^1_{\et}(X_{k_{\s}},\Z/n)\to \H^0_\et(X_{k_{\s}},\Z/n)[2]\cong\Z/n[2]$$
in $D(\Gamma, \Z/n)$ corresponding to the complex 
$\tau^{[0,1]}\RGamma_{\et}(X_{k_{\s}},\Z/n)$. The $\Z/n$-linear dual of the map $\delta_X$ shifted by $2$ is a map $\delta_X^{\vee}[2]:\Z/n\to (A[n])[2]$ which is the following composition
$$
\Z/n\xrightarrow{\Delta}\End_{\Z/n}(A[n])=\Hom_{\Z/n}(A[n],\Z/n)\otimes^L_{\Z/n}A[n]\xrightarrow{\delta_X\otimes\id}\Z/n[2]\otimes_{\Z/n}^L A[n],
$$
where $\Delta$ sends $1\in\Z/n$ to $\id_{A[n]}\in\End_{\Z/n}(A[n])$. After applying the derived functor of the functor of $\Gamma$-invariants, we see that 
$\id_{A[n]}\in\End_{\Z/n}(A[n])^\Ga$ goes to $\delta_X\in \H^2(k,A[n])$.

The map $\delta_X\otimes\id_{A[n]}$ can be identified with the connecting map of the complex $\tau^{[0,1]}\RGamma_\et(X_{k_{\s}},\Z/n)\otimes^L_{\Z/n}A[n]\simeq \tau^{[0,1]}\RGamma_\et(X_{k_{\s}},A[n])$. Thus the differential
$$d\colon\H^1_\et(X_{k_\s},A[n])^\Ga \to\H^2_\et(k,A[n])$$ of the spectral sequence
$$E_2^{p,q}=\H^p(k, \H^q_\et(X_{k_\s},A[n]))\to\H^{p+q}_\et(X,A[n])$$
is induced by $\delta_X\otimes \id_{A[n]}$. 
There are natural isomorphisms of $\Ga$-modules
$$\H^1_\et(X_{k_\s},A[n])\cong \H^1_\et(A_{k_\s},A[n])\cong\H^1_\et(A_{k_\s},\Z/n)\otimes_{\Z/n}A[n]\cong \End_{\Z/n}(A[n]).$$
In the notation of \cite[Proposition~3.2.2]{S01}, under these isomorphisms, 
$\id_{A[n]}\in \End_{\Z/n}(A[n])$ corresponds to the class
$\tau\in \H^1_\et(X_{k_\s},A[n])$ of the $X_{k_\s}$-torsor for $A[n]$ given by
(\ref{torsor}) after an identification of $X_{k_\s}$ with $A_{k_\s}$. (Note that $\tau$ does not depend on the choice of a $k_\s$-point of $X$, 
thus $\tau\in \H^1_\et(X_{k_\s},A[n])^\Ga$.) 
This gives $d(\tau)=\delta_X$. 

On the other hand, by \cite[Proposition~3.2.2]{S01}, which is a restatement
of \cite[Proposition~V.3.2.9]{G71}, we have $d(\tau)=\beta_n([X])$. \hfill $\Box$

\brem{\rm \label{section remark}
Using \cite[Lemma 2.4.5]{S01} as in \cite[Proposition 2.2]{hs}, one sees easily that 
$\beta_n([X])\in \H^2(k, A[n])$ can also be interpreted as the class of the group extension
\begin{equation}
0\to A[n]\to G_{X,n}\to \Ga\to 0 \label{GGG}
\end{equation}
obtained from the fundamental exact sequence
\begin{equation}
1\to\pi^{\et}_1(X_{k_{\s}})\to \pi^{\et}_1(X)\to \Ga\to 1\label{fundX}
\end{equation} 
by pushing out along the surjection 
$\pi_1^{\et}(X_{k_{\s}})\cong\pi_1^{\et}(A_{k_{\s}})\to A[n]$. This can be used to give
an alternative proof of Theorem \ref{1}. On the one hand, arguing as in Lemma \ref{odin}
we see that the spectral sequence (\ref{ssX}) is canonically identified
with the Hochschild--Serre spectral sequence of continuous group cohomology 
\begin{equation}
\label{HS for Y}
\H^i(k, \H^j_{\cont}(\pi_1^{\et}(X_{k_{\s}}),\Z/n))\Rightarrow \H^{i+j}_{\cont}(\pi_1^{\et}(X),\Z/n)
\end{equation}
defined by (\ref{fundX}). On the other hand, a theorem of Hochschild and Serre 
describes the differential $\delta^{0,1}_2$ of (\ref{HS for Y}) as the map sending
$x\in \Hom_\Ga(A[n],\Z/n)\cong \H_\et^1(X_{k_\s},\Z/n)^\Ga$ to the class
in $\H^2(k,\Z/n)$ given by the push-out of (\ref{GGG}) by $x$,
and similarly for $\delta^{p,1}_2$ for all $p\geq 0$, see \cite[Theorem 4]{HSerre}. }
\erem

Among the $k$-torsors $X$ for a semiabelian variety $A$ whose class in $\H^1(k,A)$ 
is not divisible by $n$, we have the following well-known examples in dimension 1.

\bexa{\rm (i) \label{Ex section 5}
Let $k$ be a field of characteristic not equal to $2$.
Let $X$ be the affine curve $x^2-ay^2=b$, where $a,b\in k^\times$. 
This is a torsor for the norm $1$ torus 
$$S=R^1_{k(\sqrt{a})/k}(\G_{m,k(\sqrt{a})}):=\Ker[R_{k(\sqrt{a})/k}(\G_{m,k(\sqrt{a})})\to
\G_{m,k}],$$
where the arrow is given by the norm ${\rm N}\colon k(\sqrt{a})\to k$. Hilbert's Theorem 90 gives an
isomorphism $\H^1(k,S)\cong k^\times/{\rm N}(k(\sqrt{a})^\times)$. 
This group is annihilated by 2, so if $[X]\neq 0$ then $[X]$ is not divisible by 2. 
The class $[X]\in \H^1(k,S)$ is zero if and only if the projective conic
$x^2-ay^2=bz^2$ has a $k$-point, that is, if and only if the symbol $(a,b)\in\H^2(k,\Z/2)$ is zero.
In fact, we have an isomorphism of group $k$-schemes $S[2]\cong\Z/2$, and 
$\delta_X=\beta_2([X])=(a,b)\in \H^2(k,\Z/2)$, as follows, for example, from 
\cite[Proposition 7.1.11]{CTS21}.

\smallskip

\noindent (ii) Let $k=\R$ and let $X$ be a smooth projective curve of genus 1 over $\R$
such that $X(\R)=\emptyset$. Let $E$ be the Jacobian of $X$. The group $\H^1(\R,E)$ is 
annihilated by 2, and since $[X]\neq 0$, we see that $\delta_X=\beta_2([X])\neq 0$.
}
\eexa

\brem{\rm
(i) More generally, if $X$ is a geometrically connected scheme over $k$, then the map $\Z/n\cong \H^0_{\et}(X_{k_{\s}},\Z/n)\to \RGamma_{\et}(X_{k_{\s}},\Z/n)$ admits a splitting in $D(\Gamma,\Z/n)$ if $X$ has a $0$-cycle of degree coprime to $n$. Indeed, for a $0$-cycle $\sum\limits_{i}a_i[Z_i]$ with each $f_i:Z_i\hookrightarrow X$ isomorphic to $\Spec(L_i)$,
where $L_i$ is a finite extension of $k$, consider the induced maps $$r_i:\RGamma_{\et}(X_{k_{\s}},\Z/n)\xrightarrow{f_i^*}\RGamma(\Spec(L_i\otimes_k k_{\s}),\Z/n)\cong \bigoplus\limits_{L_{i,\s}\hookrightarrow k_{\s}}\Z/n\xrightarrow{\Sigma}\Z/n.$$ Here the direct sum is taken over all embeddings of the maximal separable subextension $L_{i,\s}$ of $L_i$ into $k_{\s}$. The map $\Sigma$ is $\Gamma$-equivariant because the Galois group permutes the
summands in the direct sum. The map on $\H^0$ induced by $r_i$ is then the multiplication by the degree $[L_{i,\s}:k]$, so the sum $\sum a_i\cdot [L_i:L_{i,\s}]\cdot r_i$ is a map $\RGamma_{\et}(X_{k_{\s}},\Z/n)\to \Z/n$ inducing multiplication by the degree of the cycle $\sum\limits_i a_i[Z_i]$.

\noindent (ii) In particular, if $X$ is a geometrically connected smooth proper variety over $k$, then the map $\Z/n\to \RGamma_{\et}(X_{k_{\s}},\Z/n)$ admits a splitting if $n$ is coprime to the Euler characteristic $\chi(X_{k_{\s}})$ (defined, e.g., as $\sum_i (-1)^i \dim_{\Q_{\ell}}\H^i_{\et}(X_{k_{\s}},\Q_{\ell})$ for any $\ell\neq {\rm char}(k)$). Indeed, the top Chern class $c_{\dim X}(T_X)\in {\rm CH}^{\dim X}(X)$ of the tangent bundle is a $0$-cycle of degree $\chi(X)$, by the Hirzebruch--Riemann--Roch formula. By (i), such a cycle gives a map $\RGamma_{\et}(X_{k_{\s}},\Z/n)\to \Z/n$ which induces $\chi(X)\cdot\id$ on $\H^0$.

}
\erem

\section{Curves}

\newcommand{\bPic}{{\bf Pic}}
Let $C$ be a smooth, projective, geometrically integral curve over a field $k$.
Let $J:={\bf Pic}_{C/k}^0$ be the Jacobian of $C$ and let $X:={\bf Pic}_{C/k}^1$ 
be the Albanese torsor of $C$.
The canonical map $C\to X$ gives isomorphisms of $\Ga$-modules
$$\H^1_\et(C_{k_\s},\Z/n)\cong \H^1_{\et}(X_{k_{\s}},\Z/n)\cong \H^1_\et(J_{k_\s},\Z/n).$$
We use this identification in the following proposition.

\bpr \label{curves: main}
Let $n$ be an integer not divisible by ${\rm char}(k)$.

{\rm (a)} For all $i\in \Z$ and all $p\geq 0$ the differential  $\delta_{2,C}^{p,1}$ of the spectral sequence
$$E^{p,q}_2=\H^p(k,\H^q_\et(C_{k_\s},\mu_n^{\otimes i}))\Rightarrow\H^{p+q}_\et(C,\mu_n^{\otimes i}), $$
is equal to $\delta_{2,X}^{p,1}$, and is induced by the class $\beta_n([X])\in \H^2(k,J[n])$, as in Theorem~\ref{1};

{\rm (b)} the differential $\delta^{0,2}_{2,C}$ of the spectral sequence
\begin{equation}
E^{p,q}_2=\H^p(k,\H^q_\et(C_{k_\s},\mu_n))\Rightarrow\H^{p+q}_\et(C,\mu_n)
\label{ss mu}
\end{equation}
sends the generator of $\H^2(C_{k_\s},\mu_n)^\Ga\cong\H^2(C_{k_\s},\mu_n) \cong\Z/n$
to the image of $[X]$ under the map 
$$\H^1(k,J)\xrightarrow{\beta_n} \H^2(k,J[n])\cong\H^2(k,\H^1_\et(J^\vee_{k_\s},\mu_n))\cong 
\H^2(k,\H^1_\et(J_{k_\s},\mu_n)).$$
Here the last isomorphism is induced by the principal polarisation of $J$. Consequently, for all $p\geq 0$ the differential $\delta^{p,2}_{2,C}$ is given by cupping with this class in 
$\H^2(k, \H^1_{\et}(J_{k_{\s}},\mu_n))$. 
\epr

{\em Proof.} 
The canonical map $C\to X$ induces a map 
$\RGamma_{\et}(X_{k_{\s}},\mu_n^{\otimes i})\to \RGamma_{\et}(C_{k_{\s}},\mu_n^{\otimes i})$ in $D(\Gamma,\Z/n)$, which induces isomorphisms
of cohomology groups in degrees 0 and 1. This gives (a).

Poincar\'e duality states that in $D(\Gamma,\Z/n)$ the objects $\Hom_{\Z/n}(\RGamma_{\et}(C_{k_{\s}},\Z/n),\mu_n)$ and 
$\RGamma_{\et}(C_{k_{\s}},\Z/n)[2]$ are canonically isomorphic. Since each cohomology group of $\RGamma_{\et}(C_{k_{\s}},\Z/n)$ is a free $\Z/n$-module, this implies a duality on truncations: the objects
$\Hom_{\Z/n}(\tau^{[0,1]}\RGamma_{\et}(C_{k_{\s}},\Z/n),\mu_n)$ and 
$\tau^{[1,2]}\RGamma_{\et}(C_{k_{\s}},\Z/n)[2]$ are canonically isomorphic.
Therefore the differential
$\delta^{0,2}_{2,C}$ of the spectral sequence with coefficients $\Z/n$ is given by
the class of a 2-extension in $\Ext^2_k(\Z/n,\H^1_\et(J^\vee_{k_\s},\mu_n))$ 
that is dual to the extension
defining $\delta^{0,1}_{2,C}$. It follows that $\delta^{0,2}_{2,C}$ sends the generator
of $\H^2_\et(C_{k_\s},\mu_n)^\Ga\cong\Z/n$ to the image of $[X]$ in 
$\H^2(k,J[n])\cong\H^2(k,\H^1_\et(J^\vee_{k_\s},\mu_n))$. \hfill $\Box$

\bexa{\rm \label{real curves odd genus}
Let $C$ be a smooth, projective, geometrically integral curve over $\R$ such that $C(\R)=\emptyset$.
If the genus of $C$ is odd, then the Albanese torsor is non-trivial, that is, $[X]\neq 0$ in $\H^1(\R,J)$,
see \cite[Proposition 3.3 (2)]{GH81}.
Since $\H^1(\R,J)$ is annihilated by 2, we have $\beta_2([X])\neq 0$, thus for $n=2$ the differential
$\delta^{0,2}_{2,C}$ is non-zero. See Examples \ref{suslin lemma remark} and \ref{real curves even genus} below for the case of even genus.
}\eexa

\bexa{\rm \label{suslin lemma remark}
The differentials on the 3rd page of the Hochschild--Serre spectral sequence
of a smooth projective curve can be non-zero. Indeed, consider
the spectral sequence (\ref{ss mu}) for $n=2$, where $C$ is the conic $ax^2+by^2=z^2$
with $a,b\in k^\times$, ${\rm char}(k)\neq 2$.
Then $\H^1_\et(C_{k_\s},\Z/2)=0$, so all differentials on the 2nd page are zero.
On the 3rd page we have a differential
$$\delta_{3,C}^{0,2}\colon \H^2_\et(C_{k_\s},\Z/2)^\Ga\to \H^3(k,\Z/2).$$
Suslin's lemma \cite[Lemma 1]{Sus82} states that the generator of
$\H^2_\et(C_{k_\s},\Z/2)^\Ga\cong \Z/2$ goes to the symbol $(a,b,-1)$. 
In particular, if $k=\R$ and $a=b=-1$, then $\delta_{3,C}^{0,2}\neq 0$.
}\eexa

\bexa{\rm \label{general} 
More generally, let $X$ be a smooth, projective, geometrically integral variety over $\R$ such that $X(\R)=\emptyset$. It is well known that the spectral sequence 
\begin{equation}
E^{p,q}_2=\H^p(\R,\H^q_\et(X_{\C},\Z/2))\Rightarrow\H^{p+q}_\et(X,\Z/2) \label{ssC}
\end{equation}
has non-zero differentials \cite[Corollary 2.6 (b)]{AV64}, \cite[Proposition 2.2]{Cox79}. 
The following quick argument was communicated to us by Vadim Vologodsky. The Artin comparison quasi-isomorphism $\RGamma_{\et}(X_{\C},\Z/2)\simeq \RGamma_{\mathrm{sing}}(X(\C),\Z/2)$ can be made $\Gal(\C/\R)$-equivariant with the Galois group acting on the scheme $X_{\C}$ by functoriality, and on the topological space $X(\C)$ by continuous automorphisms of $\C$. Hence $\RGamma_{\et}(X,\Z/2)$ can be identified with $\RGamma(\Gal(\C/\R),\RGamma_{\mathrm{sing}}(X(\C),\Z/2))$ which is the equivariant cohomology of the topological space $X(\C)$ with respect to the action of $\Gal(\C/\R)$. 
Since $X(\R)=\emptyset$, the Galois action on $X(\C)$ is free, and hence this equivariant cohomology coincides with the cohomology of the quotient space $X(\C)/\Gal(\C/\R)$ which is necessarily concentrated in finitely many degrees. But if the spectral sequence (\ref{ssC}) had no non-zero differentials, the cohomology of the complex $\RGamma_{\et}(X,\Z/2)$ would be non-zero in infinitely many degrees.}
\eexa

One can describe the differentials $\delta_3^{p,2}$ in the spectral sequence (\ref{ss mu}) in general, assuming that all differentials on the $2$nd page vanish. It follows from Proposition \ref{curves: main} that the term $E^{0,2}_2=\H^2_{\et}(C_{k_{\s}},\mu_n)^{\Gamma}\simeq \Z/n$ remains intact on the $3$rd page if and only if the class of Albanese torsor $[\bPic^1_{C/k}]\in \H^1(k, J)$ is divisible by $n$. 
Let us assume that this is the case and choose a $k$-torsor for $J$, which we denote by 
$\bPic^{1/n}_{C/k}$, together with an isomorphism $[n]_*\bPic^{1/n}_{C/k}\simeq \bPic^1_{C/k}$. Equivalently, we have a pushout diagram of extensions of group schemes

\begin{equation}\label{curves: div by n torsor diagram}
\begin{tikzcd}
J\arrow[r] & A\arrow[r,"{\pi}"] & \Z \\
J\arrow[r]\arrow[u, equal] & \bPic_{C/k}\arrow[u,"g"]\arrow[r] & \Z\arrow[u,"{n}"]
\end{tikzcd}
\end{equation}
with $\bPic^{1/n}_{C/k}$ isomorphic to the fibre of $\pi$ over $1\in\Z$.

The Picard stack of $C$ is a $\G_m$-gerbe over the Picard scheme $\bPic_{C/k}$ and it has the structure of a group stack compatible with the group structure on $\bPic_{C/k}$. In particular, it defines a degree $2$ extension of Galois modules
$$
\delta:\bPic_{C/k}(k_{\s})\to k_{\s}^{\times}[2].
$$
By definition, the map $g$ in the diagram (\ref{curves: div by n torsor diagram}) defines a short exact sequence of Galois modules
\begin{equation}\label{curves: div by n exact sequence}
0\to \bPic_C(k_{\s})\xrightarrow{g} A(k_{\s})\to\Z/n\to 0
\end{equation}
and we define a degree $3$ extension as the composition
\begin{equation}\label{curves: deg 3 extension over Z}
\Z/n\to \bPic_C(k_{\s})[1]\xrightarrow{\delta[1]} k_{\s}^{\times}[3]
\end{equation}
with the first map being the connecting morphism in $D(\Gamma,\Z)$ induced by (\ref{curves: div by n exact sequence}). 

Recall that for an object $M\in D(\Gamma,\Z)$ there are natural isomorphisms for all $i$:
\begin{equation}\label{curves: mod n adjunction}
\H^i(k,M\otimes^{L}_{\Z}\Z/n)\cong \Ext^{i+1}_{\Gamma}(\Z/n, M)
\end{equation}
Indeed, the right-hand side can be calculated as $\H^{i+1}(k,\mathrm{RHom}_{\Z}(\Z/n,M))$ and $\mathrm{RHom}_{\Z}(\Z/n,M)\in D(\Gamma,\Z)$ is quasi-isomorphic to $M\otimes^{L}_{\Z}\Z/n[-1]$.

For a Galois module $M$, which is $n$-divisible as an abelian group, there is a natural isomorphism 
$\Ext^i_{\Gamma}(\Z/n, M)\cong \Ext^i_{\Gamma}(\Z,M[n])=\H^i(k,M[n])$. 
Let $c_{1/n}\in \H^3(k,\mu_n)$ be the image of the composition of maps 
in (\ref{curves: deg 3 extension over Z}) under this isomorphism for $M=k_{\s}^{\times}$.
As a consequence of the following proposition, the class $c_{1/n}$ does not depend on the choice of the torsor $\bPic^{1/n}_{C/k}$.
\bpr\label{curves: 3rd page main}
Let $\bPic_{C/k}^{1/n}$ be a $k$-torsor for $J$ such that 
$[n]_*[\bPic_{C/k}^{1/n}]=[\bPic_{C/k}^1]$ in $\H^1(k, J)$. Then all differentials on the second page of the Hochschild--Serre spectral sequence
$$
E_2^{p,q}=\H^p(k,\H^q_{\et}(C_{k_{\s}},\mu_n))\Rightarrow \H^{p+q}_{\et}(C,\mu_n)
$$
are zero. The only non-zero differentials on the $3$rd page
\begin{multline*}
E_{3}^{p,2}=\H^p(k, \H^2_{\et}(C_{k_{\s}},\mu_n))=\H^p(k,\Z/n)\to \\ \to E_3^{p+3,0}=\H^{p+3}(k, \H^0_{\et}(C_{k_{\s}},\mu_n))=\H^{p+3}(k,\mu_n)
\end{multline*}
are induced by the class $c_{1/n}$.
\epr
{\em Proof.} Under the assumption that $[\bPic_{C/k}^1]$ is divisible by $n$, all differentials on the second page vanish by Proposition \ref{curves: main}. To access the differentials on the next page, we relate the extensions between the cohomology modules of the complex $\RGamma_{\et}(C_{k_{\s}},\mu_n)$ to those of $\RGamma_{\et}(C_{k_{\s}},\G_m)$. 

The complex $\tau^{\leq 2}\RGamma_{\et}(C_{k_{\s}},\G_m)$
has non-zero cohomology modules only in degrees $0$ and $1$, isomorphic to $k_{\s}^{\times}$ and $\bPic_{C/k}(k_{\s})$, respectively, cf.~\cite[Theorem 5.6.1 (iv)]{CTS21}. The degree $2$ extension between them is exactly the class $\delta$. According to the Kummer sequence, applying the functor $-\otimes_{\Z}^L\Z/n:D(\Gamma,\Z)\to D(\Gamma,\Z/n)$ to the object $\RGamma_{\et}(C_{k_{\s}},\G_m)$ gives $\RGamma(C_{k_{\s}},\mu_n)[1]$, and the extension
$$
\H^0_\et(C_{k_{\s}},\G_m)\to \tau^{\leq 2}\RGamma_{\et}(C_{k_{\s}},\G_m)\to \H^1_\et(C_{k_{\s}},\G_m)[-1]
$$
is carried to
$$
\H^0_\et(C_{k_{\s}},\mu_n)[1]\to \RGamma_{\et}(C_{k_{\s}},\mu_n)[1]\to (\tau^{[1,2]}\RGamma_{\et}(C_{k_{\s}},\mu_n))[1].
$$
Hence the extension 
\begin{equation}\label{curves: ext to compute}
\tau^{[1,2]}\RGamma_{\et}(C_{k_{\s}},\mu_n)[2]\to \mu_n[3]
\end{equation} 
arising from the complex $\RGamma_{\et}(C_{k_{\s}},\mu_n)$ is the result of applying the functor 
$-\otimes_{\Z}^L\Z/n$ to the map $\delta:\bPic_{C/k}(k_{\s})\to k_{\s}^{\times}[2]$.

Under the assumption that $[\bPic_{C/k}^1]$ is divisible by $n$, the natural map 
\begin{equation}\label{curves: mod n map to split}
\tau^{[1,2]}\RGamma_{\et}(C_{k_{\s}},\mu_n)[2]\simeq \bPic_{C/k}(k_{\s})\otimes^L_{\Z}\Z/n\to \Z/n 
\end{equation} 
has a section in $D(\Gamma,\Z/n)$ and the task of computing the $3$rd differential is equivalent to computing the composition of this section with (\ref{curves: ext to compute}).

Recall from (\ref{curves: mod n adjunction}) that for any Galois module $M$ there is a natural isomorphism $\Ext^i_{\Gamma}(\Z/n, M[1])\cong \H^i(k, M\otimes^L_{\Z}\Z/n)$. The morphism $h:\Z/n\to \bPic_{C/k}(k_{\s})[1]$ corresponding to the extension $0\to \bPic_{C/k}(k_{\s})\to A(k_{\s})\to \Z/n\to 0$ is carried to a class in $\H^0(k, \bPic_{C/k}(k_{\s})\otimes_{\Z}^L\Z/n)$ which defines a splitting of (\ref{curves: mod n map to split}) because the composition $$\Z/n\xrightarrow{h}\bPic_{C/k}(k_{\s})[1]\xrightarrow{\deg[1]}\Z[1]$$ is the Bockstein map. 

The composition of this splitting with (\ref{curves: ext to compute}) is then the map $\Z/n\to \mu_n[3]$ whose class in $\H^3(k,\mu_n)$ corresponds to the composition $\Z/n\xrightarrow{h} \bPic_{C/k}(k_{\s})[1]\xrightarrow{\delta[1]}k_{\s}^{\times}[3]$, and we arrive at the definition of the class $c_{1/n}$, as desired.
\hfill$\Box$.

\medskip

If the torsion order of the Albanese torsor $[\bPic^1_{C/k}]\in \H^1(k,J)$ is coprime to $n$, then we can make the above formula for the $3$rd differential more explicit. The map $\delta$ arising from the Picard stack of $C$ induces a map $\bPic_{C/k}(k)\to \Br(k)$ sending a rational point to the obstruction to lifting it to an actual line bundle on $C$. We denote by $\Bock_{\G_m}:\Br(k)=\H^2(k,\G_m)\to \H^3(k,\mu_n)$ the connecting homomorphism of the Kummer sequence.

\bco\label{curves: prime to n 3rd page}
Suppose that $[\bPic^1_{C/k}]\in \H^1(k,J)$ is annihilated by an integer $m$ coprime to $n$. Then all differentials on the second page of the Hochschild--Serre spectral sequence
$$
\H^p(k,\H^q_{\et}(C_{k_{\s}},\mu_n))\Rightarrow \H^{p+q}_{\et}(C,\mu_n)
$$
are zero. Let $x\in\bPic_{C/k}^d(k)$, where
$d\equiv 1\bmod n$. The only non-zero differentials
$$
E_{3}^{p,2}=\H^p(k,\Z/n)\to \H^{p+3}(k,\mu_n)=E_3^{p+3,0}
$$
on the $3$rd page are given by cupping with $\Bock_{\G_m}(\delta(x))\in \H^3(k,\mu_n)$.
\eco

{\em Proof.} Since $n$ is invertible modulo $m$, the class $[\bPic^1_{C/k}]$ is divisible by $n$, so we are in the setup of Proposition \ref{curves: 3rd page main}. Specifically, choose $m'$ to be any integer  such that $m'\cdot n\equiv 1\bmod m$ and let $\bPic_{C/k}^{1/n}:=[m']_*\bPic^1_{C/k}$. The point $x$ defines a map of Galois modules $f_x:\Z\to \bPic_{C/k}(k_{\s})$ sending $1\in \Z$ to $x$. The composition
$$
\Z\xrightarrow{f_x}\bPic_{C/k}(k_{\s})\xrightarrow{\delta}k_{\s}^{\times}[2]\xrightarrow{\Bock_{\G_m}}\mu_n[3]
$$
corresponds under the isomorphism $\H^3(k,\mu_n)\simeq \Ext^3_{\Gamma}(\Z/n,k_{\s}^{\times})$ to the composition
\begin{equation}\label{curves: n divisible albanese deg 3 ext formula}
\Z/n\to\Z[1]\xrightarrow{f_x[1]}\bPic_{C/k}(k_{\s})[1]\xrightarrow{\delta[1]}k_{\s}^{\times}[3].
\end{equation}
The composition of the first two maps in
(\ref{curves: n divisible albanese deg 3 ext formula}) equals 
$$\Z/n\xrightarrow{f_x\otimes^L_{\Z}\Z/n}\bPic_C(k_{\s})\otimes^L_{\Z}\Z/n\xrightarrow{\id\otimes\Bock}\bPic_C(k_{\s})[1].$$ Since $f_x\otimes^L_{\Z}\Z/n$ is a section of the natural map $\bPic_C(k_{\s})\otimes^L_{\Z}\Z/n\to \Z/n$, the above composition (\ref{curves: n divisible albanese deg 3 ext formula}) is the map inducing the differentials on the $3$rd page of our spectral sequence, as in the proof of Proposition \ref{curves: 3rd page main}. \hfill$\Box$.

\bexa{\rm \label{real curves even genus}
Let $C$ be a smooth, projective, geometrically integral curve over $\R$ such that $C(\R)=\emptyset$.
The spectral sequence
$$E_{2}^{p,q}=\H^p(\R,\H^q_{\et}(C_{\C}, \G_{m})) \Rightarrow\H^{p+q}_{\et}(C, \G_{m})$$
gives rise to the short exact sequence
$$0\to\Pic(C)\to\Pic(C_\C)^{\Gal(\C/\R)}=\bPic_{C/\R}(\R)\xrightarrow{\delta}\Br(\R)\to 0,$$
as follows from \cite[Proposition 2.2 (2)]{GH81}.
If the genus of $C$ is even, then the Albanese torsor is trivial, that is, 
$[\bPic^1_{C/\R}]=0\in \H^1(\R,J)$,  see \cite[Proposition 3.3 (1)]{GH81}. 
Moreover, by the second statement of 
\cite[Proposition 2.2 (2)]{GH81}, for any $x\in\bPic^1_{C/\R}(\R)$ we have
$\delta(x)\neq 0$ in $\Br(\R)$. 
By Corollary \ref{curves: prime to n 3rd page} the 3rd page differentials
$\H^p(\R,\Z/2)\to \H^{p+3}(\R,\Z/2)$ are cup-products with the generator
$(-1,-1,-1)$ of $\H^3(\R,\Z/2)$, so they are isomorphisms $\Z/2\tilde\lra\Z/2$ for all $p\geq 0$.
}\eexa

\brem{\rm One can also give an explicit degree $3$ Yoneda model for the extension $\Z/n\to \mu_n[3]$ discussed above, at least after applying the forgetful functor $D(\Gamma,\Z/n)\to D(\Gamma,\Z)$.

Let $\Div_{C_{k_{\s}}}$ be the group of divisors on $C_{k_{\s}}$, that is, the free abelian group with the set of closed points of $C_{k_{\s}}$ as the basis. 
The exact sequence
\begin{equation}
0\to k_{\s}^{\times}\xrightarrow{}k_{\s}(C)^{\times}\xrightarrow{}\Div_{C_{k_{\s}}}\to \bPic_{C/k}(k_{\s})\to 0
\end{equation} 
represents the extension class of $\tau^{\leq 1}\RGamma_{\et}(C_{k_{\s}},\G_m)$,
see \cite[Proposition 5.4.5, Remark 5.4.6]{CTS21}.
From diagram (\ref{curves: div by n torsor diagram}) we see that multiplication by $n$ on 
$\bPic_{C/k}$ factors as 
$$\bPic_{C/k}\hookrightarrow A\xrightarrow{f}\bPic_{C/k}$$ for some map $f$. \newcommand{\tDiv}{\widetilde{\Div}}Let $\tDiv$ be the $\Gamma$-module defined as the pushout
$$
\begin{tikzcd}
\tDiv\arrow[r]\arrow[d] & \Div_{C_{k_{\s}}}\arrow[d] \\
A(k_{\s})\arrow[r, "f"] & \bPic_{C/k}(k_{\s})
\end{tikzcd}
$$
We will view elements of $\tDiv$ as pairs $(\alpha,D)\in A(k_{\s})\oplus \Div_{C_{k_{\s}}}$ satisfying $f(\alpha)=[D]$. We can then form the complex
\begin{equation}\label{curves: degree 3 yoneda extension}
0\to\mu_n\to k_{\s}(C)^{\times}\xrightarrow{a}k_{\s}(C)^{\times}\oplus \Div_{C_{k_{\s}}}\xrightarrow{b}\tDiv\xrightarrow{c}\Z/n\to 0.
\end{equation}
Here the map $a$ sends a rational function $\varphi$ to the pair $(\varphi^n,\div(\varphi))\in k_{\s}(C)^{\times}\oplus\Div_{C_{k_{\s}}}$. The map $b$ on $k_{\s}(C)^{\times}$ is induced from the map $-\div:k_{\s}(C)^{\times}\to \Div_{C_{k_{\s}}}$, and on $\Div_{C_{k_{\s}}}$ the map $b$ sends a divisor $D$ to the element of $\tDiv$ given by $(g([D]),n\cdot D)$. Finally, $c$ is the composition $\tDiv\to A(k_{\s})\to \mathrm{coker}(g)\simeq\Z/n$.

The sequence (\ref{curves: degree 3 yoneda extension}) is exact and hence defines a map $\Z/n\to\mu_n[3]$ in $D(\Gamma,\Z)$. Unwinding the definitions one checks that it is equal to the image of $c_{1/n}\in \H^3(k,\mu_n)$ under the map $\H^3(k,\mu_n)\to \Ext^3_{\Gamma,\Z}(\Z/n,\mu_n)$.
}\erem

We finally explain why Corollary \ref{curves: prime to n 3rd page} implies Suslin's lemma 
\cite[Lemma 1]{Sus82} discussed in Example \ref{suslin lemma remark}. 

\ble \label{last lemma}
Suppose that ${\rm char}(k)\neq 2$.
Consider the extension of $\Ga$-modules
\begin{equation}
0\to\Z/2\to\mu_4\to \Z/2\to 0. \label{mu_4}
\end{equation}
For even $n$ the connecting map $\H^n(k,\Z/2)\to\H^{n+1}(k,\Z/2)$ of $(\ref{mu_4})$
is the cup-product with the class of $-1$ in $\H^1(k,\Z/2)=k^\times/k^{\times 2}$.
\ele
{\em Proof.}
The spectral sequence $\H^p(k,\Ext_\Z^q(\Z/2,\Z/2))\Rightarrow\Ext^{p+q}_{\Ga,\Z}(\Z/2,\Z/2)$ gives an exact sequence
$$0\to \H^1(k,\Z/2)\to\Ext^1_{\Ga,\Z}(\Z/2,\Z/2)\to \Ext^1_\Z(\Z/2,\Z/2).$$
As pointed out in \cite[the proof of Lemma 4]{Sus82},
the difference between the image of $[-1]\in \H^1(k,\Z/2)$ in 
$\Ext^1_{\Gamma,\Z}(\Z/2,\Z/2)$ and the class of (\ref{mu_4}) is the class of the extension where $\Z/4$ is equipped with the trivial Galois action
\begin{equation}
0\to\Z/2\to\Z/4\to\Z/2\to 0. \label{Z/4}
\end{equation}
By the Milnor--Bloch--Kato conjecture,
the natural map
$$\H^n(k,\mu_4^{\otimes n})\to\H^n(k,\mu_2^{\otimes n})=\H^n(k,\Z/2)$$
is surjective. For even $n$ we have $\mu_4^{\otimes n}\cong \Z/4$, so 
the connecting map of (\ref{Z/4}) is zero on $\H^n(k,\Z/2)$. \hfill$\Box$

\medskip
If $C$ is a conic, then $\delta$ sends $1\in \Z\cong\Pic(C_{k_\s})^\Ga\cong\bPic_{C/k}(k)$ to
the class of the associated quaternion algebra in $\Br(k)$ (see
\cite[Proposition 7.1.3]{CTS21}), that is, to the image of the symbol $(a,b)$ under the
natural map $\H^2(k,\Z/2)\to \H^2(k,\G_m)$. 
By Lemma \ref{last lemma} the image of $(a,b)$ under the 
Bockstein map $\H^2(k,\G_m)\to\H^3(k,\Z/2)$ is $(a,b,-1)$. We of course do not need to invoke the Milnor--Bloch--Kato conjecture for this particular computation, because the Kummer classes 
of $a$ and $b$ in $\H^1(k,\Z/2)$ lift to $\H^1(k,\mu_4)$, so $(a,b)\in \H^2(k,\Z/2)$ lifts
to $\H^2(k, \mu_4^{\otimes 2})\cong\H^2(k,\Z/4)$.

\bigskip

Department of Mathematics, Massachusetts Institute of Technology, Cambridge, MA
02139, USA

\medskip

\texttt{alexander.petrov.57@gmail.com}
\bigskip

Department of Mathematics, South Kensington Campus, Imperial College London
SW7~2AZ United Kingdom \  \ and \  \ 
Institute for the Information Transmission Problems,
Russian Academy of Sciences, Moscow 127994 Russia

\medskip

\texttt{a.skorobogatov@imperial.ac.uk}

\end{document}